\documentclass{aptpub}

\RequirePackage[OT1]{fontenc}
\RequirePackage[numbers]{natbib}
\RequirePackage[colorlinks,citecolor=blue,urlcolor=blue]{hyperref}

\newtheorem{notation}{Notation}

\usepackage{graphicx}

 \usepackage{mathptmx}      
\usepackage{latexsym}

\usepackage[ansinew]{inputenc}

\usepackage{lmodern} 

\usepackage{graphicx}
\usepackage{wrapfig}
\usepackage{caption}
\usepackage{subfig}
\usepackage{wrapfig}

\usepackage{pgfplots}

\usepackage{amsfonts}
\usepackage{dsfont}

\usepackage{calc}
\usepackage{xcolor}
\usepackage{tikz}
\usetikzlibrary{decorations.markings}
\usetikzlibrary{arrows}
\tikzstyle{vertex}=[circle, draw, inner sep=0pt, minimum size=6pt]

\tikzstyle{vertexr}=[fill=red, circle, draw, inner sep=0pt, minimum size=6pt]

\tikzstyle{vertexb}=[fill=blue, circle, draw, inner sep=0pt, minimum size=6pt]

\tikzstyle{vertexg}=[fill=green, circle, draw, inner sep=0pt, minimum size=6pt]

\tikzstyle{vertexy}=[fill=yellow, circle, draw, inner sep=0pt, minimum size=6pt]

\usepackage{color}

\definecolor{gray50}{gray}{0.5}

\newcommand{\Pp}{\mathcal{P}}
\newcommand{\cible}{\overline{\vec{x}}}

\newcommand{\Pro}{\mathbb{P}}

\newcommand{\hyp}{\mathcal{H}_N}

\shorttitle{Branching random walks on binary strings}
\authornames{I. Balelli, V. Mili\v{s}i\'c and G. Wainrib}
\begin{document}

\title{Branching random walks on binary strings\\ for evolutionary processes}



\authorone[Universit\'e Paris 13,  Laboratoire Analyse G\'eom\'etrie et Applications]{Irene Balelli}
\authorone[Universit\'e Paris 13,  Laboratoire Analyse G\'eom\'etrie et Applications]{Vuk Milisic}
\authortwo[Ecole Normale Sup\'erieure, D\'epartement d'Informatique]{Gilles Wainrib}
\addressone{99, avenue Jean-Baptiste Cl\'ement
93430 - Villetaneuse - France}
\addresstwo{45 rue d'Ulm, 75005 - Paris - France.}






\begin{abstract}
In this article, we study branching random walks on graphs modeling division-mutation processes inspired by adaptive immunity. We  
apply the theory of expander graphs on  mutation rules in evolutionary processes and  obtain estimates for the cover times of the branching random walk. This analysis reveals an unexpected saturation phenomenon~:  increasing the mutation rate above a certain threshold does not enhance the speed of state-space exploration.
\end{abstract}





\section{Introduction}\label{intro}

The aim of this article is to understand interactions between mutation and division in evolutionary processes. In particular, we are interested in analyzing  characteristic time-scales for which a certain proportion of possible traits is expressed in the population~: starting from a single individual, what would be the typical time until a finite proportion of the traits are covered by the exponentially increasing population? In the models we consider, traits are represented as vertices of the $N$-dimensional hypercube, and the choice of a mutation rule corresponds to the prescription of a graph structure. The division-mutation process is then modeled as a Branching Random Walk (BRW) on this graph.
A division event is always associated to mutation, meaning that the newborn particles move neighboring nodes according to a given  mutation rule. 
We consider two kinds of branching processes~: a simple branching random walk (also called COBRA walk \cite{duttacoalescing,cooper2016coalescing}) where two or more particles having the same trait coalesce into a single one, and the branching random walk with multiplicity, where we also take into account the number of individuals sharing the same trait within the population. 
The choice of a mutation rule influences the graph structure, and we show that the theory of expander graphs leads to new results on the typical time-scales of the state-space exploration. 
To our knowledge, the link between expander graphs theory and evolutionary process is new.\\

The motivation behind this work is the study of Antibody Affinity Maturation (AAM). This is a key process of the adaptive immune system, which allows to create specific high-affinity antibodies against pathogens that threaten a given organism. 
Antibodies are proteins which are secreted by B-cells, special lymphocytes  trained to recognize the presented antigen \cite{LSom}. 
This process takes place in Germinal Centers (GCs) \cite{KMMurPTraMWal,victora2014snapshot}, 
where activated B-cells proliferate, mutate and differentiate. B-cells recognize the antigen thanks to transmembrane proteins called B-cell Receptors (BCRs) \cite{de2015dynamics,kringelum2013structural}. The mutational mechanism that B-cells undergo during a GC reaction is called Somatic Hypermutation (SHM)~: it affects, at a very high rate, the DNA encoding for the specific portion of the BCR involved in the binding with the antigen, called Variable (V) region \cite{teng2007immunoglobulin}. \\

\begin{figure}[ht!]
\subfloat[ ]{\includegraphics[width = 1.5in]{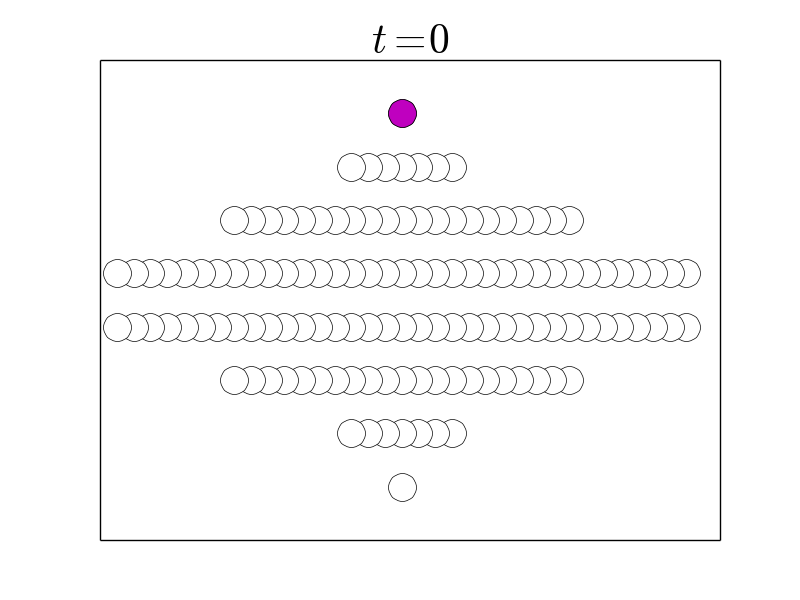}\label{fig1:intro}}
       ~ 
\subfloat[ ]{\includegraphics[width = 1.5in]{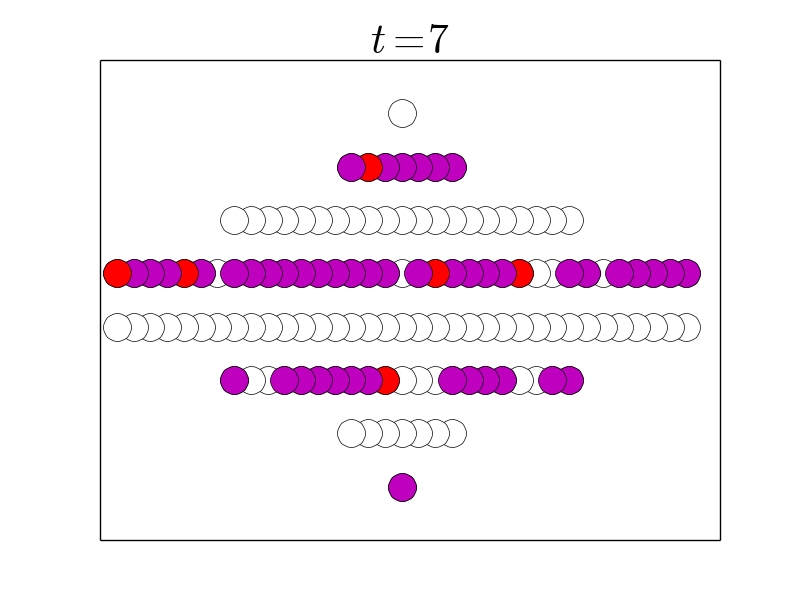}\label{fig2:intro}}
       ~ 
\subfloat[ ]{\includegraphics[width = 1.5in]{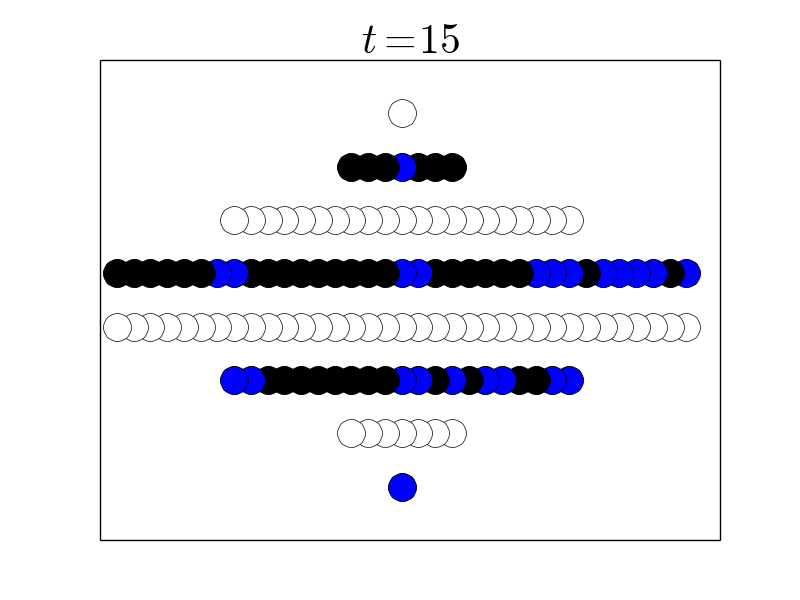}\label{fig2:intro}}
\caption{Simulation of the exploration of the state-space of all possible traits, considering division and simple switch-type mutations, \emph{i.e.} a mutation consists in the switch of a randomly chosen bit. Traits are represented by 7-length binary strings and we give a picture of the evolution of the process at the beginning (a) after 7 time steps (b) and after 15 time steps (c). The process starts with a single individual whose trait corresponds to $\vec x_0=[0,0,0,0,0,0,0]$, represented by the top circle. Starting from the top, on each line we arrange all nodes having an increasing distance from $\vec x_0$, so that, for instance, the bottom circle corresponds to the node $[1,1,1,1,1,1,1]$. Different colors denote a different number of individuals sharing the same trait. In particular in magenta we plotted nodes with at most 4 individuals lying on them, in red between 5 and 9, in blue between 300 and 499 and in dark more than 500 individuals. }\label{fig:intro}
\end{figure}

A certain number of mathematical models and results about GC reaction and AAM already exists. In particular, T. B. Kepler and A. S. Perelson in \cite{kepler1993cyclic,kepler1993somatic} proposed deterministic population dynamics models for SHM and AAM, 
considering for the first time the hypothesis of the existence of a recycling mechanisms of B-cells during GC reaction. 
This mechanism has now been confirmed by experiments \cite{victora2010germinal}. In \cite{oprea1997somatic,ASPerGWeis,forrest2014population,DIbePKMai} the authors introduced and discussed several immunological problems, such as the size of the repertoire, or the strength of antigen-antibody binding, while providing as well suitable mathematical tools. More recently, other articles have focused on biologically detailed models of the GC reaction (\emph{e.g.} \cite{meyer2002mathematical}), in particular with an agent-based modeling framework (\cite{meyer2012theory}, mostly analyzed through extensive numerical simulations).
In 2015 the journal \textit{Philosophical Transactions of the Royal Society B} has entirely dedicated an issue to the dynamics of antibody repertoires. For instance, in \cite{elhanati2015inferring,mccoy2015quantifying,cobey2015evolution} the authors developed and applied modern statistical methods to investigate selection on BCRs and infer B-cell sequence evolution.
We are interested in studying from an analytical point of view evolutionary pathways of BCRs  during SHM. Here and in \cite{balelli2015branching,iba.3}, we
provide some significant building blocks in this direction and study their mathematical features.\\

Besides the biological motivations, the class of models studied in this article is interesting from a mathematical point of view, as it is a discrete-time BRW on graphs, a type of branching process which has not been deeply investigated so far to our knowledge, despite its growing number of applications. Since the first articles about branching processes in the 50's and 60's \cite{kendall1948generalized,bellman1952age,ikeda1968branching1,ikeda1968branching2,ikeda1969branching3}, this class of stochastic processes has been used in various situations to model biological, genetic, physical, chemical or technological processes. For example branching processes can model the dynamics of population in genetics \cite{sawyer1976branching}, or the spread of a piece of data, a rumor or a virus \cite{ball1995strong}. Most of the works that have been published so far are not interested in studying these processes on graphs. Nevertheless, in some recent papers \cite{bertacchi2006weak,bertacchi2009characterization} the authors considers BRWs on multigraphs and mostly focus on weak and strong survival conditions. Dutta C. \emph{et al} in \cite{duttacoalescing} exhibits bounds on  cover times 
for COBRA walks on trees, grids, and expander graphs (useful later in our analysis) in the context of gossip propagation. Results on expander graphs have been  improved in \cite{cooper2016coalescing}  using a new duality relation between the COBRA walk and a discrete epidemic process. 
Another field of recent interest is the study of BRWs in random environnement. 
We refer, for example, to \cite{albeverio2000annealed}, where the authors study local and total particle populations or to \cite{machado2003branching} where conditions for recurrence and transience (almost surely wrt the random environment) are found, for the discrete-time BRW on a rooted tree with random environment. Branching annihilating RWs have  been extensively studied in last years due to their applications in biological, chemical, physical and economical systems \cite{marro2005nonequilibrium,cardy1996theory,cardy1998field}. In \cite{szabo2000branching} the authors consider these processes on random regular graphs, which they study using Monte Carlo simulations and the generalized mean-field analysis.\\

In this paper, we focus on BRWs on $\{0,1\}^N$ with constant division rate 2 (except for Section \ref{sec5:1}), inspired by cellular division. The coupling of branching mechanism and random walk necessarily implies an important speedup in the characteristic time-scales for the exploration of the state-space. Typically, for the simple random walk on the $N$-dimensional hypercube, the addition of a branching process enables a speedup from a time $\mathcal O(2^N)$ to $\mathcal O(N)$ (Section \ref{sec3:2}). Of course this has a cost~: considering a branching process means also to produce new individuals at each time step. Indeed, in a time $T=\mathcal O(N)$ we have $2^T$ individuals (in the case in which multiplicity is taken into account; $\leq 2^N$ otherwise), as we do not consider here neither selection nor death. The mutation rule, which defines the structure of the graph, also determines the ability of the BRW in covering the vertices of the graph. In particular, using expansion properties, in Section \ref{sec3} we prove that the best result we can obtain in a time $\mathcal O(N)$ for finite connected expander graphs over the state-space $\{0,1\}^N$, is to cover a half of the graph. \\

Moreover, our mathematical analysis of the cover times has revealed an interesting phenomenon concerning the impact of the mutation rate on the exploration speed. Intuitively, one would suggest that increasing the number of mutations at each division  would result in a BRW with a faster exploration time-scale. However, we show in Section \ref{sec3:3} the existence of an early saturation phenomenon : when increasing from one to two mutations, the exploration is indeed faster, but allowing more than two mutations (up to $N$) modifies only marginally the exploration speed.\\

In Section \ref{sec:1} we state the main definitions and notations setting up a general mathematical framework. 
Section \ref{sec:bip} contains preliminary results concerning generic BRWs on graphs and their possible bipartite structure. 
Bipartiteness actually influences  the dynamic of the branching process. In Section \ref{sec3}, we establish quantitative results concerning the portion of the state-space invaded in $\mathcal O(N)$ for two different kinds of BRWs  (Theorems \ref{lempartial} and \ref{lempartialbis}). In order to do so, we need to determine some characteristics of the graphs, in particular their expansion properties. 
These results provide quantitative estimates of the typical time-scale for state-space exploration resulting from the interaction between division and mutation. 
Then, in Section \ref{sec5}, we propose some extensions of the model. 
In particular, we introduce the BRW with multiplicity and obtain the transition matrix related to 
the  number of individuals carrying a given trait together with their limiting distribution. 
We investigate as well how this distribution can change by introducing a division rate, and provide comparisons between different mutation/division models.
In this way,  theoretical results presented in previous sections are displayed in a wider context.
Finally, in Section \ref{sec:conclusion} 
we conclude with a brief summary of this work and discuss the biological setting in which it is embedded 
justifying our hypotheses. We present as well consequences of our results and 
 discuss possible  improvements in order to  cover the state-space faster in time, or to drive the covering to main interest areas of the graph.

\section{Definitions and Notations}\label{sec:1}

We start this section with some definitions and notations, establishing an elementary mathematical framework for the modeling of antibody affinity maturation in the germinal center. \\

We first assume that it is possible to classify the amino acids, which determine the che\-mi\-cal properties of both epitope and paratope, into 2 classes, typically positively charged and negatively charged. Henceforth BCRs and antigen are represented by binary strings of a same length $N$, hence, the state-space of all possible BCR configurations is $\{0,1\}^N$ (we refer to \cite{balelli2015branching} for more details). 

\begin{definition}
We denote by $\hyp$ the standard $N$-dimensional hypercube. BCR and antigen configurations are represented by vertices of $\hyp$, denoted by $\vec x_i$ with $1\leq i\leq2^N$, or sometimes simply by their indices.
\end{definition}

In this paper we introduce and discuss models including mutation and division. Mathematically, this gives rise to BRWs on $\{0,1\}^N$. The structure of the graph depends then on the mutation rule we consider.\\ 

We suppose that there is a single B-cell entering the GC reaction. At each time step, each B-cell divides and mutates according to a given mutational rule. A mutation corresponds to a jump on a neighbor node.

\begin{definition}\label{defn:neighgeneral}
Given $\vec{x}_i$, $\vec{x}_j\, \in \{0,1\}^N$, we say that $\vec{x}_i$ and $\vec{x}_j$ are neighbors, and denote $\vec{x}_i \sim \vec{x}_j$, if there exists at least one edge (or loop) between them.
\end{definition}

We are mostly interested in studying the variation of the number of expressed traits within the population, as a result of the interaction between division and mutation. 
In this paper we refer to two different kinds of BRWs: the simple $c$-BRW (also called coalescing BRW \cite{duttacoalescing}) and the $c$-BRW with multiplicity.

\begin{definition}[Simple $\mathbf c$-BRW]\label{def:simpleBRW}
The process starts at an arbitrary node (representing the BCR of a B-cell entering the process of division and mutation during the GC reaction), labelled as active. If at time $t$ node $\vec x_i$ is active (\emph{i.e.} the trait $\vec x_i$ is expressed in the GC population at time $t$), then at time $t+1$ it chooses $c$ of its neighbors, independently and with replacement, to become active, while $\vec x_i$ becomes inactive again (unless, of course, another active node at time $t$ chooses it). In this model, the number of times a node is chosen to be active is not taken into account. We suppose $c>1$, otherwise the BRW simply becomes a RW.
\end{definition}

\begin{definition}[$\mathbf c$-BRW with multiplicity]\label{def:BRWmult}The process starts with a B-cell entering the process of mutation and division, lying on an arbitrary node which corresponds to its trait. At each time step a particle lying on a certain node $\vec x_i$ of $\{0,1\}^N$ gives rise to $c$ daughter cells, with $c>1$, and die. Each one of the $c$ newborn particles choses a neighbor node, independently and with replacement, and move on it. More than one particle can lie on the same vertex of $\hyp$, and each one divides at each time step.
\end{definition}

\begin{notation}\label{noataion1}
Let $S\,\subseteq\,V$ be a subset of vertices of a graph $G=(V,E)$. Then we denote by $\mathcal{N}(S)$ the set of the neighbors of all vertices in $S$. We denote by $|S|$ and $|\mathcal{N}(S)|$ the number of vertices in $S$ and in $\mathcal{N}(S)$ respectively. $\mathcal{N}(S)$ may include also some vertices in $S$.
\end{notation}

\begin{notation}
Given a simple $c$-BRW on a generic graph $G$, for all $t\geq0$ we note by $S_t$ the set of all active nodes at time $t$ and by $\mathcal{N}(S_t)$ the set of all the neighbors of the vertex set $S_t$.
\end{notation}

The structure of the graph and consequently the dynamics of the BRW on it depend on the introduced mutation rule, which is defined thanks to the transition probability matrix.

\begin{definition}
Let $\mathcal M$ be the transition probability matrix of a graph $G$. We denote the BRW referring to $\mathcal M$ and with constant division rate $c$ by $c$-BRW-$\mathcal M$.
\end{definition}

In particular, we refer to two mutational rules (see \cite{balelli2015branching} for more details). Here below we give the definitions of the corresponding transition probability matrices.

\begin{definition}\label{def:Pp}
For all $\vec{x}_i$, $\vec{x}_j\,\in\,\hyp$~:
\begin{equation*}
\Pro(\vec{X}_n=\vec{x}_j\,|\,\vec{X}_{n-1}=\vec{x}_i)=:p(\vec{x}_i,\vec{x}_j)=\left\{\begin{array}{ll}
1/N & \textrm{if}\quad \vec{x}_j \sim \vec{x}_i \\
0 & \textrm{otherwise}
\end{array}\right.
\end{equation*}
\end{definition}
Matrix $\Pp:=\left(p(\vec{x}_i,\vec{x}_j)\right)_{\vec{x}_i,\vec{x}_j\in\hyp}$ gives to $\{0,1\}^N$ the structure of a standard $N$-dimensional hypercube.\\

We further introduce another transition matrix, which models a mutation rule in which up to $k$ symbols of the string are independently mutated at each division:
\begin{definition}\label{def:Ppn}
Let $k\;\in\;\{1,\dots, N\}$, $\Pp^{(k)}:=\displaystyle\frac{1}{k}\sum_{i=1}^k \Pp^i$, $\Pp$ given by Definition \ref{def:Pp}.
\end{definition}

We finally recall the definition of Hamming distance, which measures, in our model, the affinity between two traits  \cite{balelli2015branching}:

\begin{definition}\label{def:hamm}
For all $\vec x =(x_1,\dots,x_N)$, $\vec y=(y_1,\dots,y_N)\,\in\,\{0,1\}^N$, their Hamming distance is given by:
\begin{displaymath}
h(\vec{x},\vec{y})=\sum_{i=1}^N\delta_i\qquad\textrm{where}\qquad \delta_i=\left\{\begin{array}{ll}
1 & \textrm{if}\quad x_i \neq y_i \\
0 & \textrm{otherwise}
\end{array}\right.
\end{displaymath}
\end{definition}

\begin{definition}\label{def:affhamming}
For all $\vec x_i\,\in\,\{0,1\}^N$, its affinity with a given vertex $\cible$, $\textrm{aff}(\vec x_i, \cible)$ is  given by $\textrm{aff}(\vec x_i, \cible):=N-h(\vec x_i,\cible)$, where $h(\cdot,\cdot):(\{0,1\}^N\times\{0,1\}^N)\rightarrow\{0,\dots, N\}$ returns the Hamming distance.
\end{definition}

\section{$\mathbf c$-BRW on graphs and bipartiteness}\label{sec:bip}

The bipartiteness deeply influences the characteristics of the BRW and, in particular, its possibility of covering all nodes of the graph simultaneously at a certain time. 

\begin{definition}
A graph $G=(V,E)$ is bipartite if there exists a partition of the vertex set $V=V_1\sqcup V_2$, s.t. every edge connects a vertex in $V_1$ to a vertex in $V_2$.
\end{definition}
We emphasize the relations between a generic $c$-BRW on a given graph $G=(V,E)$, with $c\geq2$, and the eventual bipartite structure of the above-mentioned graph.

\subsection{$c$-BRW on bipartite graphs}

Let us consider a simple $c$-BRW on a generic bipartite graph $G_b(V_1\sqcup V_2,E)$. Instead of considering a single random active node at the beginning, we suppose that the process starts with a given initial distribution $\mathbf{p}$ of the active set. The results presented in this section do not change if we consider a $c$-BRW with multiplicity instead of a simple $c$-BRW. The fact that the trials are made with replacement does not have any consequences either.

\begin{proposition}
If the initial distribution $\mathbf{p}$ is concentrated on $V_1$ or on $V_2$ then $|S_t|\leq\max_{i=1,2}{(|V_i|)}$ for all $t\geq0$, otherwise $S_t=V_1\sqcup V_2$ for some $t>0$ with positive probability.
\end{proposition}

\proof
The proof is a direct consequence of the bipartite structure of $G_b$. Let us suppose, without loss of generality, that $\mathbf{p}$ is concentrated on $V_1$. Then, due to the bipartite structure of $G_b$ after an even number of steps we have necessarily $S_{2t}\subseteq V_1$, while after an odd number of steps we have $S_{2t+1}\subseteq V_2$, and so the first statement is proven.\\
If, on the contrary, $\mathbf{p}$ is not concentrated on $V_1$ nor on $V_2$, then for all $t\geq0$ we have a positive probability that $S_t=S_{t,1}\sqcup S_{t,2}$ with $S_{t,1}\subseteq V_1$ and $S_{t,2}\subseteq V_2$, and consequently, w.p.p. we have $S_t=V_1\sqcup V_2$ for some $t>0$.
\endproof

\begin{remark}
This qualitative result does not change if we take into account the number of times a node is chosen to become active for the next time step or if we decide to make trials without replacement: these choices only have effects on the speed of the covering.
\end{remark}

\subsection{$c$-BRW on non-bipartite connected graphs}\label{ss:BRW_nonbip}

Let us now consider a non-bipartite connected graph $G=(V,E)$. We recall a classical result about bipartite graphs \cite{salvatore2007bipartite}, which will be useful later:

\begin{proposition}\label{nonbip0}
A graph is bipartite if and only if it has no odd cycles.
\end{proposition}

We shall prove the following statement:

\begin{theorem}\label{nonbip}
Given a $c$-BRW on a finite non-bipartite connected graph, then, w.p.p., there exists a time $t>0$ such that $S_t=V$.
\end{theorem}

The proof of this theorem is based on the following three lemmas:

\begin{lemma}\label{nonbip1}
If $G=(V,E)$ is a finite connected graph, then, independently from the initial distribution, $\forall\;\vec x_i\,\in\,V$ there exists a time $t<\infty$ s.t. $\vec x_i\,\in\,S_t$.
\end{lemma}

In other words, if the graph is finite and connected, then each node will be activated by the BRW at least once in a finite time interval.

\proof
The hitting time of the $c$-BRW to reach any node of a finite connected graph, starting from every possible initial distribution is finite, thanks to the connectivity of the graph and the fact that it has a finite set of nodes. (Note that this is still true if $c=1$, \emph{i.e.} for a SRW on a finite connected graph).
\endproof

\begin{lemma}\label{nonbip2}
If there exists a time $t\geq0$ such that $\exists\;\vec x_1$, $\vec x_2$, $\vec x_1\sim \vec x_2$ and $\{\vec x_1,\vec x_2\}\,\in\,S_t$, then w.p.p. there exists a time $T>t$ s.t. $S_T=V$ (independently from the initial distribution).
\end{lemma}

This means that if at a given time $t$ we have two neighbor nodes both active, then we have a positive probability to reach $S_T=V$ later.

\proof
Let us suppose that $\vec x_1$, $\vec x_2$ are two neighbor nodes and $S_{t^\ast}=\{\vec x_1,\vec x_2\}$ (we suppose that all other nodes are non-active). Then we are able to show that w.p.p., for all $t\geq t^\ast$, $S_t\subset \mathcal{N}(S_t)$ and $S_t=\mathcal{N}(S_t)\Leftrightarrow S_t=V$, where we recall that $\mathcal{N}(S_t)$ is the set of all neighbors of $S_t$. This implies that w.p.p. the active set can always grow until we reach $S_t=V$. This result is quite intuitive, indeed if $\vec x_1\sim \vec x_2$ and $S_{t^\ast}=\{\vec x_1,\vec x_2\}$, then necessarily $S_{t^\ast}\subset \mathcal{N}(S_{t^\ast})$ and, consequently, there is a positive probability that $S_{t^\ast}\subset S_{t^\ast+1}$. That means that w.p.p. $S_{t^\ast+1}$ contains $x_1$, $x_2$ and at most $c-1$ distinct neighbors of $\vec x_1$ and $c-1$ distinct neighbors of $\vec x_2$. Then we can repeat the same argument with all the couples of neighbors active at time $t^\ast+1$ (w.p.p. all nodes in $S_{t^\ast+1}$ are neighbors two by two). Thanks to the connectivity of the graph and the fact that it is a finite graph, w.p.p. we can go on with this procedure until we reach $S_t=V$.
\endproof

\begin{lemma}\label{nonbip3}
If there exists at least an odd cycle on $G=(V,E)$, then, independently from the initial distribution, w.p.p. for a time $t\geq0$ there exist two nodes $\vec x_1$, $\vec x_2$, $\vec x_1\sim \vec x_2$ and $\vec x_1$, $\vec x_2\,\in\, S_t$.
\end{lemma}

\proof
Let us suppose that in graph $G$ there exists an odd cycle of length $2n+1$~: $C=(\vec x_{1}, \vec x_{2},\dots,\vec x_{(2n+1)},\vec x_{1})$. Lemma \ref{nonbip1} implies the existence of a time $T<\infty$ s.t. $\vec x_{1}\,\in\,S_T$. Then w.p.p. we have that $\{\vec x_{2},\vec x_{(2n-1)}\}\subseteq S_{T+1}$ (we recall that $c\geq2$). We make another step and w.p.p. we have that $\{\vec x_{1},\vec x_{3},\vec x_{2n}\}\subseteq S_{T+2}$. After $n$ steps, w.p.p. we have $\{\vec x_{(n+1)},\vec x_{(n+2)}\}\subseteq S_{T+n}$, and $\vec x_{(n+1)}\sim \vec x_{(n+2)}$, which proves Lemma \ref{nonbip3}.
\endproof

\proof[Theorem \ref{nonbip}]
Let $G=(V,E)$ be a finite non-bipartite connected graph and $C=(\vec x_{1}, \vec x_{2},\dots,\vec x_{(2n+1)},\vec x_{1})$ an odd cycle of $G$ (it necessarily exists as $G$ is non-bipartite (Proposition \ref{nonbip0})). Lemma \ref{nonbip1} assures that there exists a finite time $t_1$ s.t. $\vec x_{1}\in S_{t_1}$. Then, thanks to Lemma \ref{nonbip3}, w.p.p. there exists a time $t_2>t_1$ s.t. $\{\vec x_{(n+1)},\vec x_{(n+2)}\}\subseteq S_{t_2}$, and $\vec x_{(n+1)}\sim \vec x_{(n+2)}$. The proof of Theorem \ref{nonbip} can now be achieved by applying Lemma \ref{nonbip2}.
\endproof

\section{Portion of $\mathbf{\hyp}$ covered in $\mathbf{\mathcal O(N)}$ for the simple 2-BRW-$\mathbf{\Pp}$ and the simple 2-BRW-$\mathbf{\Pp}^{(\mathbf k)}$}\label{sec3}

In this Section our aim is to estimate the size of the active node set in a time of the order of $N$. It clearly depends on the mutational model allowed on the state-space. We can interpret it as the number of possible BCR configurations expressed in our population after $\mathcal O(N)$ mutation steps. In Table \ref{tab:pct} we summarize the main results of the current section. In Sections \ref{sec3:2} and \ref{sec3:3} we estimate the size of the active set in $\mathcal O(N)$ for the simple 2-BRW referring to $\Pp$ and $\Pp^{(k)}$ (Definitions \ref{def:Pp} and \ref{def:Ppn}). We prove that the 2-BRW-$\Pp$ covers a small portion of $\hyp$, while a half of the state-space will be covered if we take into account $\Pp^{(k)}$ as transition probability matrix, at least for $N$ big enough. \\

\begin{table}[ht!]
\caption{Summary of the main results of Sections \ref{sec3:2} and \ref{sec3:3}.}\label{tab:pct}
\centering
\begin{tabular}{| c | c |}
\hline
\textbf{Model} & \textbf{$|\mathbf{S_T}|$ in $\mathbf T=\mathbf{\mathcal O(N)}$}  \\
\hline
$\mathbf{\Pp}$ &  $|S_T|\geq2^{N-r}$, $r>\displaystyle\frac{N^2e^{-2}+N-2}{Ne^{-2}+N-2}$ \\
\hline
$\mathbf{\Pp^{(k)}}$ &  $|S_T|\geq\delta2^{N}$, $\delta\leq1/2$\\
\hline
\end{tabular}
\end{table}

In order to estimate these quantities, we apply a method used in \cite{duttacoalescing} to determine the partial cover time for expander graphs. The partial cover time corresponds to the expected time required to visit at least a certain portion of the state-space. We need to evaluate the expansion properties of the graphs described by $\Pp$ and $\Pp^{(k)}$ respectively. For this reason in Section \ref{sec3:1} we recall some definitions and results about expander graphs. For a more complete overview about this subject see \cite{hoory2006expander}. 

\subsection{Expander graphs}\label{sec3:1}

Informally, an expander graph is a graph $G=(V,E)$ which has strong connectivity properties (quantified using vertex, edge or spectral expansion).  We give some mathematical characterization of this property.\\

Unless stated otherwise, throughout this section a graph $G=(V,E)$ is a connected undirected $d$-regular graph with $|V|=n$.
\begin{definition}[$(\alpha,\delta)$-expander graph]\label{defalphadelta}
$G$ is said to be an $(\alpha,\delta)$-expander graph, with $\delta\leq1/2$, if: $\forall\,S\,\subseteq\,V\;\textrm{s.t.}\;|S|\leq\delta n\;\Rightarrow\; |\mathcal{N}(S)|\geq\alpha|S|$.
\end{definition}

In other words, an $(\alpha,\delta)$-expander graph is a graph where the set of all neighbors of each subset $S$ with at most $\delta n$ nodes, has at least $\alpha|S|$ vertices.

\subsubsection{Spectrum and expansion}

Let us denote by $A_G$ the adjacency matrix of $G$ and by $\Pp_G$ its transition probability matrix. As $G$ is a $d$-regular graph, then $\Pp_G=\frac{1}{d}A_G$. We denote by $d=\lambda_1^A\geq\lambda_2^A\geq\cdots\geq\lambda_n^A$ and $1=\lambda_1\geq\lambda_2\geq\cdots\geq\lambda_n$ the eigenvalues of $A_G$ and $\Pp_G$ respectively.

\begin{definition}
We say that $G$ is a $\lambda$ eigenvalue expander, with $\lambda<d$, if $\lambda_2^A\leq\lambda$. It is a $\lambda$ absolute eigenvalue expander if $|\lambda_2^A|$, $|\lambda_n^A|\leq\lambda$.
\end{definition}

\begin{remark}
All $d$-regular connected graphs are $\lambda_2^A$ eigenvalue expanders. Indeed under these hypotheses, the first largest eigenvalue of the adjacency matrix corresponds to $d$ and $d>\lambda_2^A\geq\lambda_i^A$ for all $i\geq2$.
\end{remark}

Then we have the following known result (first proved by R. M. Tanner in \cite{tanner1984explicit}):

\begin{theorem}[Vertex expansion]\label{thm:expansion}
Let $G$ be a $\lambda$ eigenvalue expander. Let $S\subseteq V$ s.t. $|S|\leq n/2$. Then $\mathcal{N}(S)$ is large, in particular:
$$|\mathcal{N}(S)|\geq\frac{|S|}{\frac{\lambda^2}{d^2}+\left(1-\frac{\lambda^2}{d^2}\right)\frac{|S|}{n}}$$ 
\end{theorem}

\begin{remark}
One easily notices that $\left(\frac{\lambda^2}{d^2}+\left(1-\frac{\lambda^2}{d^2}\right)\frac{|S|}{n}\right)^{-1}\to 1$ for $\lambda\to d$, 
is decreasing wrt $\lambda$. 
\end{remark}

We also give another characterization of $d$-regular expander graphs with respect to their eigenvalues.

\begin{definition}\label{def:epsexp}
We say that $G$ is an \mbox{$\varepsilon$-expander} graph, with $\varepsilon<1$, if the eigenvalues of its adjacency matrix are such that $|\lambda_i^A|\leq\varepsilon d$ for $i\geq2$.
\end{definition}

Then in particular, we have the following proposition.

\begin{proposition}\label{propmu2}
Let $G$ be not bipartite. Then $G$ is a $\lambda_2$-expander graph.
\end{proposition}

\proof
As $\Pp_G=\frac{1}{d}A_G$, we have that: $|\lambda_i^A|=|\lambda_i|\cdot d\leq\lambda_2\cdot d$, $\forall\,i\geq2$.
This is not true for bipartite graphs as their spectrum is symmetric with respect to zero. Therefore $|\lambda_n^A|=d>\lambda_2\cdot d$.
\endproof

As an immediate consequence of Theorem \ref{thm:expansion} applied to $\varepsilon$-expander graphs, we have:

\begin{proposition}\label{thm:epsexpa}
Let $G$ be a $\varepsilon$-expander graph. Then for all $S\subseteq V$ s.t. $|S|\leq\delta n$, $\delta\leq1/2$, we have:
$$|\mathcal{N}(S)|\geq\frac{|S|}{\varepsilon^2(1-\delta)+\delta}$$
\end{proposition}

Finally, let us underline the clear relation existing between Definitions \ref{def:epsexp} and \ref{defalphadelta} of $\varepsilon$-expander graphs and $(\alpha,\delta)$-expander graphs respectively:
\begin{corollary}\label{cor:epsalphadelta}
Let $G$ be not bipartite with second largest eigenvalue $\lambda_2$, $\delta\leq1/2$. Then $G$ is a $(\alpha,\delta)$-expander graph with:
$$\alpha=\frac{1}{\lambda_2^2(1-\delta)+\delta}$$
\end{corollary}

\proof
First of all, Proposition \ref{propmu2} tells us that $G$ is a $\lambda_2$-expander graph. Then the condition on $\alpha$ is given by Proposition \ref{thm:epsexpa}.
\endproof

\subsection{Simple 2-BRW-$\Pp$}\label{sec3:2}

A simple 2-BRW-$\Pp$ on $\hyp$ is a generalization of a Simple RW on $\hyp$ \cite{balelli2015branching}. We want to estimate the size of the active set in $\mathcal O(N)$ using $\Pp$ as transition probability matrix. In order to do so, we use an application of a more general method used in \cite{duttacoalescing} to evaluate partial cover times. We show that the partial cover time for the simple 2-BRW-$\Pp$ is linear in $N$, while we already know that for the SRW on $\hyp$ it is exponential in $N$ \cite{avin2004efficient}. This highlights how the branching process gives an important speedup in exploring the hypercube. This speedup in covering is not without a cost. Indeed, for a time $t$ large enough, the size of the population will be of the order of the maximal possible size of $S_t$, which is $2^{N-1}$ in this case (as $\hyp$ is bipartite) and $2^N$ in the case of the simple 2-BRW-$\Pp^{(k)}$. \\

Let us start with a preliminary result about the standard $N$-dimensional hypercube, $\hyp$.

\begin{proposition}\label{expander}
For any $N\geq1$, $\hyp$ is a $N$-regular $(r,2^{-r})$-expander graph, where $r\,\in\,\{1,\dots,N\}$, \emph{i.e.}:
$$\forall\,r\,\in\,\{1,\dots,N\},\;\forall\,S\,\subset\,\{0,1\}^N\;\textrm{s.t.}\;|S|\leq2^{N-r}\quad\Rightarrow\quad |\mathcal{N}(S)|\geq r|S|$$
\end{proposition}

Before giving the proof of Proposition \ref{expander}, let us observe the maximal number of common neighbors among two or more nodes in $\hyp$.

\begin{remark}\label{shareneig}
Two distinct vertices $\vec x_{1}$, $\vec x_{2}\,\in\,\hyp$ cannot share more than two common neighbors. More generally,  $s$ distinct vertices in $\hyp$, $\{\vec x_i\}_{1\leq i\leq s\leq2^N}$  cannot share more than $s$ common neighbors.
\end{remark}

Let $A_N$ be the standard representation of the transition probability matrix of $\hyp$, obtained recursively as follows \cite{SFFlo}:
\begin{equation*}
A_1=\left(\begin{array}{cc}
0 & 1 \\
1 & 0
\end{array}\right);\;A_N=\left(\begin{array}{c|c}
A_{N-1} & \mathcal I_{2^{N-1}} \\
\hline
\mathcal I_{2^{N-1}} & A_{N-1}
\end{array}\right), \textrm{where $\mathcal I_{2^{N-1}}$ is the $2^{N-1}$-identity matrix}
\end{equation*}
Then the result is obvious since the main diagonal of $A_{N-1}$ is composed by zeros and that $A_{N-1}$ is a symmetric matrix.

\proof[Proposition \ref{expander}]
We prove Proposition \ref{expander} by double induction on $N$ and on $r$. \\

First of all, the statement is true for $N=1$ and $r=1$, and for $N=2$ and $r\,\in\,\{1,2\}$. 
We suppose the statement is true up to dimension $N-1$ and for all $r\,\in\,\{1,\dots,N-1\}$, and we prove it for dimension $N$ and for all  $r\,\in\,\{1,\dots,N\}$.\\

If $r=N$ it is true, as $\hyp$ is a $N$-regular graph. \\
Let $r=N-1$. Then we want to show that $\forall\,S\,\subset\,\{0,1\}^N\;\textrm{s.t.}\;|S|\leq2\;\Rightarrow\;|\mathcal{N}(S)|\geq (N-1)|S|$.
If $|S|=1$, for the $N$-regularity we have necessarily: $|\mathcal{N}(S)|=N>N-1$.\\
We suppose $|S|=2$, and we consider the graph of $A_N$. \\
If we choose both vertices $\vec x_i$ with \mbox{$0\leq i \leq 2^{N-1}$}, then we know, for the induction hypothesis on $N$ and observing that the top right block of $A_N$ is an identity matrix, that~: $|\mathcal{N}(S)|\geq (N-2)|S|+|S|=(N-1)|S|$. \\
Now, let us consider two vertices $\vec x_i$ and $\vec x_j$ s.t. \mbox{$i\,\in\,\{1,\dots,2^{N-1}\}$} and  $j\,\in\,\{2^{N-1}+1,\dots,2^N\}$. If we do not want to increase considerably $|\mathcal{N}(S)|$, once $\vec x_i$ is fixed, we need to choose $\vec x_j$ so that $\vec x_i$ and $\vec x_j$ share two common neighbors (Remark \ref{shareneig}).\\

Then, at least we have $|\mathcal{N}(S)|\geq2N-2=(N-1)2=(N-1)|S|$.\\
We suppose that the statement is true for dimension $N$ and for all \mbox{$r\,\in\,\{t+1,\dots,N\}$}. We prove that it's also true for $r=t$, \emph{i.e.}:
$$\forall\,S\,\subset\,\{0,1\}^N\;\textrm{s.t.}\;|S|\leq2^{N-t}\;\Rightarrow\;|\mathcal{N}(S)|\geq t|S|$$
If $|S|\leq2^{N-(t+1)}<2^{N-t}$ then, for the induction hypothesis on $r$, we have:
$$|\mathcal{N}(S)|\geq(t+1)|S|>t|S|$$
Let us suppose $2^{N-(t+1)}<|S|\leq2^{N-t}$.\\
Again, if we choose all vertices $\vec x_i$ so that $i\,\in\,\{1,\dots,2^{N-1}\}$, for the induction hypothesis on $N$ and as $r<N-1$, we have: $|\mathcal{N}(S)|\geq (t-1)|S|+|S|=t|S|$.\\
Then, we take $S=\{\vec x_i\}_{1\leq i\leq2^{N-t}}$ so that:
$$S=S_1\sqcup S_2\quad S_1=\{\vec x_{i_1}\}_{1\leq i_1\leq2^{N-1}}\quad\textrm{and}\quad S_2=\{\vec x_{i_2}\}_{2^{N-1}+1\leq i_2\leq2^N}$$
Furthermore, we suppose: $|S_1|\leq2^{N-(t+1)}$ and $|S_2|\leq2^{N-(t+1)}$, as the other cases are less favorable, if our purpose is to minimize $|\mathcal{N}(S)|$.\\
Then we have, from the induction hypothesis on $N$, together with Remark \ref{shareneig}:
$$|\mathcal{N}(S)| = |\mathcal{N}(S_1\sqcup S_2)|\geq t|S_1| + |S_1|+ t|S_2| + |S_2| - |S| = t|S|.$$
\endproof

\begin{remark}
Considering a simple $c$-BRW-$\Pp$ starting from a single node, we have that $\mathcal{N}(S_t)\cap S_t=\emptyset$ because of the bipartite structure of the graph. This is not true for generic not-bipartite graphs (see Section \ref{sec:bip}).
\end{remark}

We start by demonstrating the following lemma:

\begin{lemma}\label{lemnu}
Given a simple 2-BRW-$\Pp$~:
$$\forall\;t\geq0\textrm{ s.t. }|S_t|\leq2^{N-r}\quad\Rightarrow\quad\E[|S_{t+1}|]\geq(1+\nu)|S_t|$$
for some constant $\nu>0$ and for $r>\frac{N^2e^{-2}+N-2}{Ne^{-2}+N-2}$.
\end{lemma}

Before demonstrating Lemma \ref{lemnu}, we prove an elementary result, that we will need later:

\begin{lemma}\label{leme}
Let $c>0$ and $a$, $b>1$ such that $a\leq b$. Then:
$$e^{-ca}+e^{-cb}<e^{-c(a-1)}+e^{-c(b+1)}$$
\end{lemma}

\proof
\begin{eqnarray*}
e^{-ca}+e^{-cb}-\left(e^{-c(a-1)}+e^{-c(b+1)}\right) & = & e^{-ca}\left(1-e^{c}\right)+e^{-c(b+1)}\left(e^{c}-1\right) \\
& = & \left(1-e^{c}\right)\left(e^{-ca}-e^{-c(b+1)}\right) \\
& < & 0
\end{eqnarray*}
since $c>0$ and $a<b+1$
\endproof

\proof[Lemma \ref{lemnu}]
Let $t\geq0$ so that $|S_t|\leq2^{N-r}$, for a certain $r\,\in\,\{1,\dots,N\}$ that we will discuss later. The claim is proved if we show:
\begin{equation}\label{eqlem1}
\E[|\mathcal{N}(S_t)-S_{t+1}|]\leq|\mathcal{N}(S_t)|-(1+\nu)|S_t|
\end{equation}
For all vertices $v\,\in\,\mathcal{N}(S_t)$, let $X_v$ be the indicator variable:
\begin{equation*}
X_v=\left\{\begin{array}{ll}
1 & \textrm{if }v\,\notin\,S_{t+1} \\
0 & \textrm{otherwise}
\end{array}\right.
\end{equation*}
Then we have: $\Pro[X_v=1]=\left(1-\frac{1}{N}\right)^{2d_v}=:p$, where $d_v$ represents the number of edges connecting $v$ to $S_t$ ($1\leq d_v\leq N$).\\
Clearly $\E[X_v]=p$. Now we have:
\begin{eqnarray*}
\E[|\mathcal{N}(S_t)-S_{t+1}|] & \leq & \E\left[\sum_{v\in \mathcal{N}(S_t)}X_v\right] = \sum_{v\in \mathcal{N}(S_t)}\left(1-\frac{1}{N}\right)^{2d_v} \\ \\
& \leq & \sum_{v\in \mathcal{N}(S_t)} e^{-\frac{2d_v}{N}}
\end{eqnarray*}
Thanks to Lemma \ref{leme}, we can claim that this expression is maximized if for any $v$ (except possibly for one) $d_v$ is either 1 or $N$. In particular let us suppose that all $d_v$ are equal to 1 or to $N$ and let us denote:
$$R_1=|\{v\,\in\,\mathcal{N}(S_t)\,|\,d_v=1\}|\quad\textrm{and}\quad R_N=|\{v\,\in\,\mathcal{N}(S_t)\,|\,d_v=N\}|$$
If we are able to demonstrate the result in this particular case, then it will be true for all possible distributions of $d_v$ in $\mathcal{N}(S_t)$.
Observing that $\sum_{v\in \mathcal{N}(S_t)}d_v=N|S_t|$ thanks to the $N$ regularity, we have:
\begin{equation*}
\left\{\begin{array}{l}
R_1+R_N=|\mathcal{N}(S_t)| \\
R_1+NR_N=N|S_t|
\end{array}\right.\Rightarrow
\left\{\begin{array}{l}
R_1 = \displaystyle\frac{N}{N-1}(|\mathcal{N}(S_t)|-|S_t|) \\ \\
R_N=\displaystyle\frac{1}{N-1}(N|S_t|-|\mathcal{N}(S_t)|) 
\end{array}\right.
\end{equation*}
Then we have:
\begin{eqnarray*}
\E[|\mathcal{N}(S_t)-S_{t+1}|] & \leq & R_1e^{-\frac{2}{N}} + R_Ne^{-2} \\ 
& = & \frac{N}{N-1}(|\mathcal{N}(S_t)|-|S_t|)e^{-\frac{2}{N}} + \frac{1}{N-1}(N|S_t|-|\mathcal{N}(S_t)|)e^{-2} 
\end{eqnarray*}
In order to obtain  \eqref{eqlem1}, we have to impose that~:
\begin{equation*}
\frac{N}{N-1}(|\mathcal{N}(S_t)|-|S_t|)e^{-\frac{2}{N}} + \frac{1}{N-1}(N|S_t|-|\mathcal{N}(S_t)|)e^{-2} \leq |\mathcal{N}(S_t)|-(1+\nu)|S_t|
\end{equation*}
This is equivalent to:
\begin{equation*}
|\mathcal{N}(S_t)|\left(1-\frac{N}{N-1}e^{-\frac{2}{N}}+\frac{1}{N-1}e^{-2}\right)+|S_t|\left(\frac{N}{N-1}e^{-\frac{2}{N}}-\frac{N}{N-1}e^{-2}-1\right)\geq\nu|S_t|
\end{equation*}
By hypothesis $|S_t|\leq2^{N-r}$, which implies $|\mathcal{N}(S_t)|\geq r|S_t|$ (Proposition \ref{expander}). Since $1-\frac{N}{N-1}e^{-\frac{2}{N}}+\frac{1}{N-1}e^{-2}>0$, the last inequality will be true if:
\begin{equation}\label{eq:passaggio}
r\left(1-\frac{N}{N-1}e^{-\frac{2}{N}}+\frac{1}{N-1}e^{-2}\right)+\left(\frac{N}{N-1}e^{-\frac{2}{N}}-\frac{N}{N-1}e^{-2}-1\right)\geq\nu
\end{equation}
And so, our aim is to find $r(N)$ s.t. for all $r>r(N)$~:
\begin{equation}\label{eqtemp}
r\left(1-\frac{N}{N-1}e^{-\frac{2}{N}}+\frac{1}{N-1}e^{-2}\right)+\left(\frac{N}{N-1}e^{-\frac{2}{N}}-\frac{N}{N-1}e^{-2}-1\right)>0
\end{equation}
And this is true iff:
\begin{equation}\label{eq:Rn}
r>\frac{Ne^{-2}+N-1-Ne^{-\frac{2}{N}}}{e^{-2}+N-1-Ne^{-\frac{2}{N}}}=:r(N)
\end{equation}
We  rearrange  \eqref{eqtemp} writting~:
$$(r-1)\left(1-\frac{N}{N-1}e^{-\frac{2}{N}}\right)-\frac{N-r}{N-1}e^{-2}>0$$
Then, since $e^{-\frac{2}{N}}\leq1-\frac{2}{N}+\frac{2}{N^2}$ (thanks to the second-order Taylor expansion with integral rest), we obtain that \eqref{eqtemp} is satisfied if:
$$(r-1)\left(1-\frac{N}{N-1}\left(1-\frac{2}{N}+\frac{2}{N^2}\right)\right)-\frac{N-r}{N-1}e^{-2}>0$$
And finally:
$$r>\frac{N^2e^{-2}+N-2}{Ne^{-2}+N-2}$$
\endproof

\begin{remark}
\begin{itemize}
\item If $N\geq2$, the condition on $r$ that we found in Lemma \ref{lemnu} is met if:
$$r>1+Ne^{-2}\left(\frac{N-1}{N-2}\right)$$
\item If $N\geq3$, then this condition is satisfied if $r>1+2Ne^{-2}$.
\end{itemize}
\end{remark}

We could also express $r$ as a function of $\nu$ (we refer to  \eqref{eq:passaggio}):\\

\begin{corollary}
$\E[|S_{t+1}|]\geq(1+\nu)|S_t|$ for some constant $\nu>0$ and for 
$$N\geq r\geq\frac{\nu(N-1)+Ne^{-2}-Ne^{-\frac{2}{N}}+N-1}{e^{-2}-Ne^{-\frac{2}{N}}+N-1}:=r_N(\nu)$$
Therefore  $|S_t|$ has an exponential growth with rate $\nu$ until it reaches the size of $2^{N-r}$ and for $r\geq r_N(\nu)$. Moreover, as expected, if we define $\nu^{\ast}$ as the bigger admissible $\nu$, \emph{i.e.} $\nu^{\ast}=\sup\{\nu\;|\;r_N(\nu)\leq N\}$, then $\nu^{\ast}\leq1$.
\end{corollary}

\proof
The proof consists in elementary computations, starting from  \eqref{eq:passaggio}. 
In particular as far as the second statement is concerned, we impose $r_N(\nu)\leq N$, and clearly this condition is satisfied iff:
$$\nu\leq N-1-Ne^{-\frac{2}{N}}$$
Then, as $e^{-\frac{2}{N}}\geq1-\frac{2}{N}$, we can conclude.
\endproof

We are now able to state the following result: 

\begin{theorem}\label{lempartial}
Given a simple 2-BRW-$\Pp$, there exists a time $T$ such that $T=\mathcal O(N)$ and with high probability $|S_T|\geq2^{N-r}$, $r$ satisfying the hypothesis of Lemma \ref{lemnu}.
\end{theorem}

\proof
The proof is a direct application of a result obtained for generic expander graphs in \cite{duttacoalescing}, Section 4. 
This result applies to our specific case thanks to Lemma \ref{lemnu}. The main idea to prove Theorem \ref{lempartial} is to describe the change in the number of active nodes as a Markov process which lower bounds the growth of the size of the active set $|S_t|$. The statement is proven for this Markov process and, consequently, it is true also for our BRW.
\endproof

\subsection{Simple 2-BRW-$\Pp^{(k)}$}\label{sec3:3}

Let us start by examining an analog of Lemma \ref{lemnu} for the 2-BRW-$\Pp^{(k)}$, where we recall that $\Pp^{(k)}=\displaystyle\frac{1}{k}\sum_{i=1}^k\Pp^i$ (Definition \ref{def:Ppn}). We show that in this case the BRW covers a significantly bigger proportion of vertices in a time $\mathcal O(N)$. We follow again the method used in \cite{duttacoalescing}.\\

\begin{figure}[ht!]
  \centering
 \includegraphics[scale=0.3]{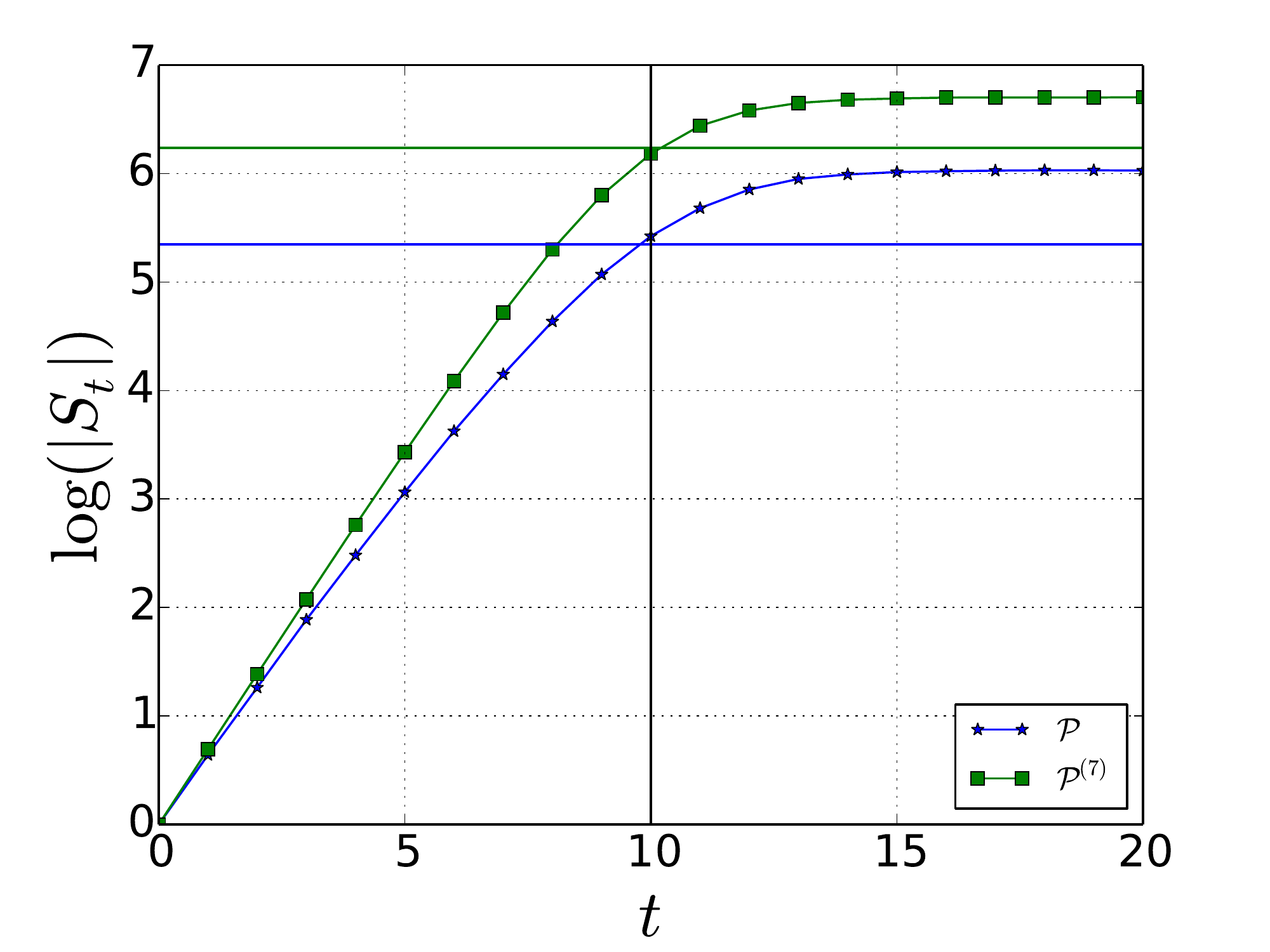}
  \caption{Evolution of the size of the active set in logarithmic scale comparing the 2-BRW for $\Pp$ (blue stars) and $\Pp^{(7)}$ (green squares) for $N=10$ (average values obtained over 40 simulations). These simulations show that the BRW referring to $\Pp^{(k)}$ can explore  the whole set of hypercube's vertices simultaneously (the graph underlying by $\Pp^{(k)}$ is not bipartite). Moreover, it is also the faster one in covering. The dark vertical line corresponds to $t=N$ and the blue and green horizontal lines represent the theoretical percentage of nodes we're supposed to cover in a time $\mathcal O(N)$, as proven in Theorem \ref{lempartial} for $\Pp$ and \ref{lempartialbis} for $\Pp^{(k)}$, for $\Pp$ and $\Pp^{(7)}$ respectively.}
  \label{fig6}
\end{figure}

First of all, we prove that the 2-BRW-$\Pp^{(k)}$ allows, for $k\geq2$, an exponential growth until it covers at least a half of the vertex set of the hypercube:

\begin{theorem}\label{lemnubis}
Given a 2-BRW-$\Pp^{(k)}$~:
$$\forall\;t\geq0\textrm{ s.t. }|S_t|\leq\delta2^N\quad\Rightarrow\quad\E[|S_{t+1}|]\geq(1+\nu)|S_t|$$
for some constant $\nu>0$, $\delta\leq1/2$ and for $N$ big enough.
\end{theorem}

In order to prove Theorem \ref{lemnubis}, we need two preliminary results (Propositions \ref{prop:Mni} and \ref{prop:Dnk}) about the characteristics of $\Pp^{(k)}$.

\begin{proposition}\label{prop:Mni}
Let $M_{N,i}:=\displaystyle\min_{j, l}{\left\{\left(\Pp^i\right)_{j,l}\,\Big|\,\left(\Pp^i\right)_{j,l}\neq0\right\}}$. We have:
$$\forall\,N\geq1,\,\forall\, i\,\in\,\{0,\dots, N\},\,M_{N,i}=i!\cdot N^{-i}$$
\end{proposition}

\proof
Due to the regularity of $\hyp$, we have that $\Pp=\frac{1}{N}A_N$. Then, $M_{N,i}=N^{-i}\min_{j\, l}$ $\left\{\left(A_N^i\right)_{j,l}\,\Big|\,\left(A_N^i\right)_{j,l}\neq0\right\}$ for all $1\leq i\leq N$, while $\Pp^0=A_N^0=I_{2^{N}}$ for all $N$ (and consequently $M_{N,0}=1$). We have now to prove that $\min_{j, l}{\left\{\left(A_N^i\right)_{j,l}\,\Big|\,\left(A_N^i\right)_{j,l}\neq0\right\}}=i!$. The proof comes directly by applying the well known result \cite{cvetkovic1995spectra}~: if $A$ is the adjacency matrix of a graph $G$, $i$ a positive integer, then the $(j,l)^{\textrm{th}}$ entry of $A^i$ corresponds to the number of $i$-length walks between vertex $j$ and $l$ in $G$.\\

First,  in the case of the $N$-dimensional hypercube, the Hamming distance between $j$ and $l$ corresponds to the length of the minimal path (\emph{i.e.} walk without loops) connecting these vertices. Moreover, because of the bipartite structure of $\hyp$ with a $i$-length walk we can not pass from $j$ to $l$ so that $h(j,l)=i-(2t+1)$, $t\geq 0$ (\emph{i.e.} with a $i$-length walk we can connect nodes having distance $k\leq i$, $k$ with the same parity as $i$). It is also clear that if $h(j,l)>i$, then there does not exists any $i$-length walk from $j$ to $l$. The minimal number of $i$-length walks to connect two nodes $j$, $l$ s.t. $h(j,l)=i-2t$, $t\geq0$ corresponds to the case $t=0$. First, if $h(j,l)=i$ we are counting the number of paths between $j$ and $l$, and that corresponds to $i!$ (we have just to choose the order of switching of the $i$ different bits). We briefly prove by combinatory arguments that given $j_1$, $l_1$, $j_2$, $l_2$ s.t. $h(j_1,l_1)=i$ and $h(j_2,l_2)=i-2$, $i\leq N$, then $\left(A^i\right)_{j_1,l_1}\leq\left(A^i\right)_{j_2,l_2}$ \emph{i.e.} $\left(A^i\right)_{j_2,l_2}\geq i!$. In order to cover a distance $i-2$ with an $i$-length walk we need to change the $i-2$ different bits in $i$ steps. Then the number of possible $i$-length walks to go from $j_2$ to $l_2$ is given by the sum for $k=0$ to $i-2$ of those walks given by the compositions of:
\begin{itemize}
\item a $k$-length path from $j_2$ to $j_{2,1}$ s.t. $h(j_{2,1},l_2)=i-2-k$~: $\binom{i-2}{k}k!$ possible choices;
\item a step from $j_{2,1}$ to $j_{2,2}$ s.t. $h(j_{2,2},l_2)=i-1-k$~: $(N-(i-2-k))$ possible choices;
\item a $(i-k-1)$-length path from $j_{2,2}$ to $l_2$~: $(i-k-1)!$ possible choices.
\end{itemize}
Finally:
\begin{equation}\label{eq:i-2}
\left(A^i\right)_{j_2,l_2}=\sum_{k=0}^{i-2}\binom{i-2}{k}k!(N-(i-2-k))(i-k-1)!
\end{equation}
We have now to prove that \eqref{eq:i-2}$\geq i!$~:
\begin{equation*}
\sum_{k=0}^{i-2}\binom{i-2}{k}k!(N-(i-2-k))(i-k-1)!=(i-2)!\sum_{k=0}^{i-2}(N-(i-2-k))(i-k-1)
\end{equation*}
And then  \eqref{eq:i-2}$\geq i!$ $\Leftrightarrow$ $\sum_{k=0}^{i-2}(N-(i-2-k))(i-k-1)\geq i(i-1)$. One can prove by an elementary computation that $\sum_{k=0}^{i-2}(N-(i-2-k))(i-k-1)=\frac{1}{6}i(i-1)(3N-2i+4)$. Consequently the result is proven if $3N-2i+4\geq6$~:
$$3N-2i+4\geq N+4 \textrm{ as $i\leq N$, and } N+4 \geq 6 \textrm{ as $N\geq2$.}$$
\endproof

Then, we give recursively the number of neighbors of each node within our graph:

\begin{proposition}\label{prop:Dnk}
Let $d_N^{(k)}$ be the number of neighbors of a generic node $l$ (including possibly $l$) in the graph corresponding to $\Pp^{(k)}$~: $d_N^{(k)}=\left|\left\{l\,\big|\,\left(\Pp^{(k)}\right)_{j,l}\neq0\right\}\right|$ for all $l\,\in\,\{1,\dots,2^N\}$ fixed. Then, $\forall\,N\geq2$~:
\begin{displaymath}
\left\{\begin{array}{l}
d_N^{(1)}=N \\ \\
d_N^{(2)}=N+d_{N-1}^{(2)} \\ \\
d_N^{(k)}=d_{N-1}^{(k-1)}+d_{N-1}^{(k)}\textrm{ for }3\leq k\leq N-1 \\ \\
d_N^{(N)}=2^N
\end{array}\right.
\end{displaymath} 
\end{proposition}

\proof
For $k=1$ and $k=N$ the proof is straightforward: if $k=1$ then we are considering the standard $N$-dimensional hypercube, and $d_N^{(1)}$ corresponds to the regularity of the graph, while if $k=N$, as we allow all possible switch-type mutations, each vertex is connected to itself and any other node within the graph. Then, as we have exactly $2^N$ vertices, $d_N^{(N)}=2^N$. In order to prove both cases $k=2$ and $3\leq k\leq N-1$ we rewrite $d_N^{(k)}$ by using powers of $A_N$. Indeed, as $\Pp^{(k)}=\frac{1}{k}\sum_{i=1}^k\left(\frac{1}{N}A_N\right)^i$, we have: $d_N^{(k)}=\left|\left\{l\,\big|\,\left(\sum_{i=1}^k A_N^i\right)_{j,l}\neq0\right\}\right|$. Proposition \ref{prop:Dnk} can now be proven by using the recursive construction of the adjacency matrix of $\hyp$ \cite{SFFlo}.
\endproof

\proof[Theorem \ref{lemnubis}]
Let $t\geq0$ so that $|S_t|\leq\delta2^{N}$, for $\delta\leq1/2$ still unknown. As we did while proving Lemma \ref{lemnu}, our aim is to show:
\begin{equation}\label{eqlembis}
\E[|\mathcal{N}(S_t)-S_{t+1}|]\leq|\mathcal{N}(S_t)|-(1+\nu)|S_t|
\end{equation}
For all vertices $v\,\in\,\mathcal{N}(S_t)$, let $X_v$ be the indicator variable:
\begin{equation*}
X_v=\left\{\begin{array}{ll}
1 & \textrm{if }v\,\notin\,S_{t+1} \\
0 & \textrm{otherwise}
\end{array}\right.
\end{equation*}
Then we have: $\Pro[X_v=1]=\displaystyle\prod_{j\sim v,\,j\in S_t}\left(1-\Pp_{jv}^{(k)}\right)^{2}=:p$. We can maximize $p$ as follows:
\begin{equation*}
p\leq \displaystyle\prod_{j\sim v,\,j\in S_t} \left(1-\frac{1}{k}\sum_{i=1}^k M_{N,i} \right)^{2}=\left(1-\frac{1}{k}\sum_{i=1}^k M_{N,i}\right)^{2d_v},
\end{equation*}
where $d_v$ represents the number of neighbors that $v$ has in $S_t$ ($1\leq d_v\leq d_N^{(k)}$).
As $\E[X_v]=p$, we have:
\begin{equation}\label{eq:pct21}
\E[|\mathcal{N}(S_t)-S_{t+1}|] \leq \sum_{v\in \mathcal{N}(S_t)}\left(1-\frac{1}{k}\sum_{i=1}^k M_{N,i}\right)^{2d_v}
\end{equation}
Denoting by $\Delta:=(1/k)\sum_{i=1}^k M_{N,i}$, we finally obtain:
\begin{equation}
\eqref{eq:pct21}\leq\sum_{v\in \mathcal{N}(S_t)} e^{-2\Delta\cdot d_v}
\end{equation}
Applying Lemma \ref{leme} this expression is maximized if for any $v$ (except possibly for one) $d_v=1$ or $d_v=d_N^{(k)}$. In particular let us suppose that all $d_v$ are equal to 1 or to $d_N^{(k)}$ and let us denote $R_1=|\{v\,\in\,\mathcal{N}(S_t)\,|\,d_v=1\}|$ and $R_2=|\{v\,\in\,\mathcal{N}(S_t)\,|\,d_v=d_N^{(k)}\}|$. We demonstrate the statement in this particular case. As $\sum_{v\in \mathcal{N}(S_t)}d_v=d_N^{(k)}|S_t|$~:
\begin{equation*}
\left\{\begin{array}{l}
R_1+R_2=|\mathcal{N}(S_t)| \\ \\
R_1+d_N^{(k)}R_2=d_N^{(k)}|S_t|
\end{array}\right.\Rightarrow
\left\{\begin{array}{l}
R_1 = \displaystyle\frac{d_N^{(k)}}{d_N^{(k)}-1}(|\mathcal{N}(S_t)|-|S_t|) \\ \\
R_N=\displaystyle\frac{1}{d_N^{(k)}-1}(d_N^{(k)}|S_t|-|\mathcal{N}(S_t)|) 
\end{array}\right.
\end{equation*}
Then we have:
\begin{equation*}
\E[|\mathcal{N}(S_t)-S_{t+1}|] \leq \frac{d_N^{(k)}}{d_N^{(k)}-1}(|\mathcal{N}(S_t)|-|S_t|)e^{-2\Delta} + \frac{1}{d_N^{(k)}-1}(N|S_t|-|\mathcal{N}(S_t)|)e^{-2\Delta d_N^{(k)}} 
\end{equation*}
Equation \eqref{eqlembis} will be satisfied if:
\begin{equation*}
\frac{d_N^{(k)}}{d_N^{(k)}-1}(|\mathcal{N}(S_t)|-|S_t|)e^{-2\Delta} + \frac{1}{d_N^{(k)}-1}(N|S_t|-|\mathcal{N}(S_t)|)e^{-2\Delta d_N^{(k)}} \leq |\mathcal{N}(S_t)|-(1+\nu)|S_t|
\end{equation*}
As the graph we are considering is a $\lambda_{N,2}^{(k)}$-expander graph (where $\lambda_{N,2}^{(k)}$ $=\frac{N-2}{2k}$ $\left(1-\left(\frac{N-2}{N}\right)^k\right)$ is the second largest eigenvalue of $\Pp_N^{(k)}$ \cite{balelli2015branching}), and applying Proposition \ref{thm:epsexpa}, the last inequality will be true if:
\begin{equation}\label{eqtempbis}
\alpha_N^{(k)}\left(1-\frac{d_N^{(k)}\cdot e^{-2\Delta}}{d_N^{(k)}-1}+\frac{e^{-2\Delta d_N^{(k)}}}{d_N^{(k)}-1}\right)+\left(\frac{d_N^{(k)}\cdot e^{-2\Delta}}{d_N^{(k)}-1}-\frac{d_N^{(k)}\cdot e^{-2\Delta d_N^{(k)}}}{d_N^{(k)}-1}-1\right)>0,
\end{equation}
where $\alpha_N^{(k)}=\displaystyle\frac{1}{\delta\left(1-{\lambda_{N,2}^{(k)}}^2\right)+{\lambda_{N,2}^{(k)}}^2}$. That means
\begin{equation}\label{eq:delnk}
\delta<\frac{e^{-2\Delta d_N^{(k)}}-d_N^{(k)}e^{-2\Delta}+d_N^{(k)}-1}{\left(1-{\lambda_{N,2}^{(k)}}^2\right)\left(d_N^{(k)}e^{-2\Delta d_N^{(k)}}-d_N^{(k)}e^{-2\Delta}+d_N^{(k)}-1\right)}-\frac{{\lambda_{N,2}^{(k)}}^2}{1-{\lambda_{N,2}^{(k)}}^2}:=\delta_N^{(k)}
\end{equation}
Finally, let us prove that for fixed $k\geq2$, $\delta_N^{(k)}$ tends to 1 for $N$ going to infinity. Indeed we have:
\begin{itemize}
\item Let $k\geq2$~: $\Delta=\frac{1}{k}\sum_{i=1}^kM_{N,i}=\frac{1}{k}\left(\frac{1}{N}+\frac{2}{N^2}\right)+\frac{1}{k}\sum_{i=1}^k\frac{i!}{N^i}$. And then, for $N\to\infty$, $\Delta\sim \mathcal O\left(\frac{1}{N}\right)$
\item For fixed $k$, $d_N(k)$ is monotonically increasing:
\begin{itemize}
\item $k=1\Rightarrow d_N^{(1)}=N>N-1=d_{N-1}^{(1)}$;
\item $k=2\Rightarrow d_N^{(2)}=N+d_{N-1}^{(2)}>d_{N-1}^{(2)}$;
\item $3\leq k\leq N-1\Rightarrow d_N^{(k)}=d_{N-1}^{(k-1)}+d_{N-1}^{(k)}>d_{N-1}^{(k)}$;
\end{itemize}
\item Let $k\geq2$~: $d_{N}^{(k)}\geq d_N^{(2)}$. By definition: $d_N^{(2)}=N+d_{N-1}^{(2)}=\sum_{i=0}^{N-3}(N-i)+d_2^{(2)}=\frac{N^2+N+2}{2}$. Therefore, for fixed $k\geq2$, $\Delta d_N^{(k)}$ tends to infinity for $N\to \infty$.
\end{itemize}
Finally we have, for $k\geq 2$ fixed:
$$\delta_N^{(k)}=\frac{e^{-2\Delta d_N^{(k)}}-d_N^{(k)}\left(e^{-2\Delta}-1\right)-1}{\left(1-{\lambda_{N,2}^{(k)}}^2\right)\left(d_N^{(k)}e^{-2\Delta d_N^{(k)}}-d_N^{(k)}\left(e^{-2\Delta}-1\right)-1\right)}-\frac{{\lambda_{N,2}^{(k)}}^2}{1-{\lambda_{N,2}^{(k)}}^2}\to1\quad\textrm{for $N\to\infty$}$$

\begin{figure}[ht!]
  \centering
  \subfloat[$N=7$]{\label{fig6a}\includegraphics[scale=0.29]{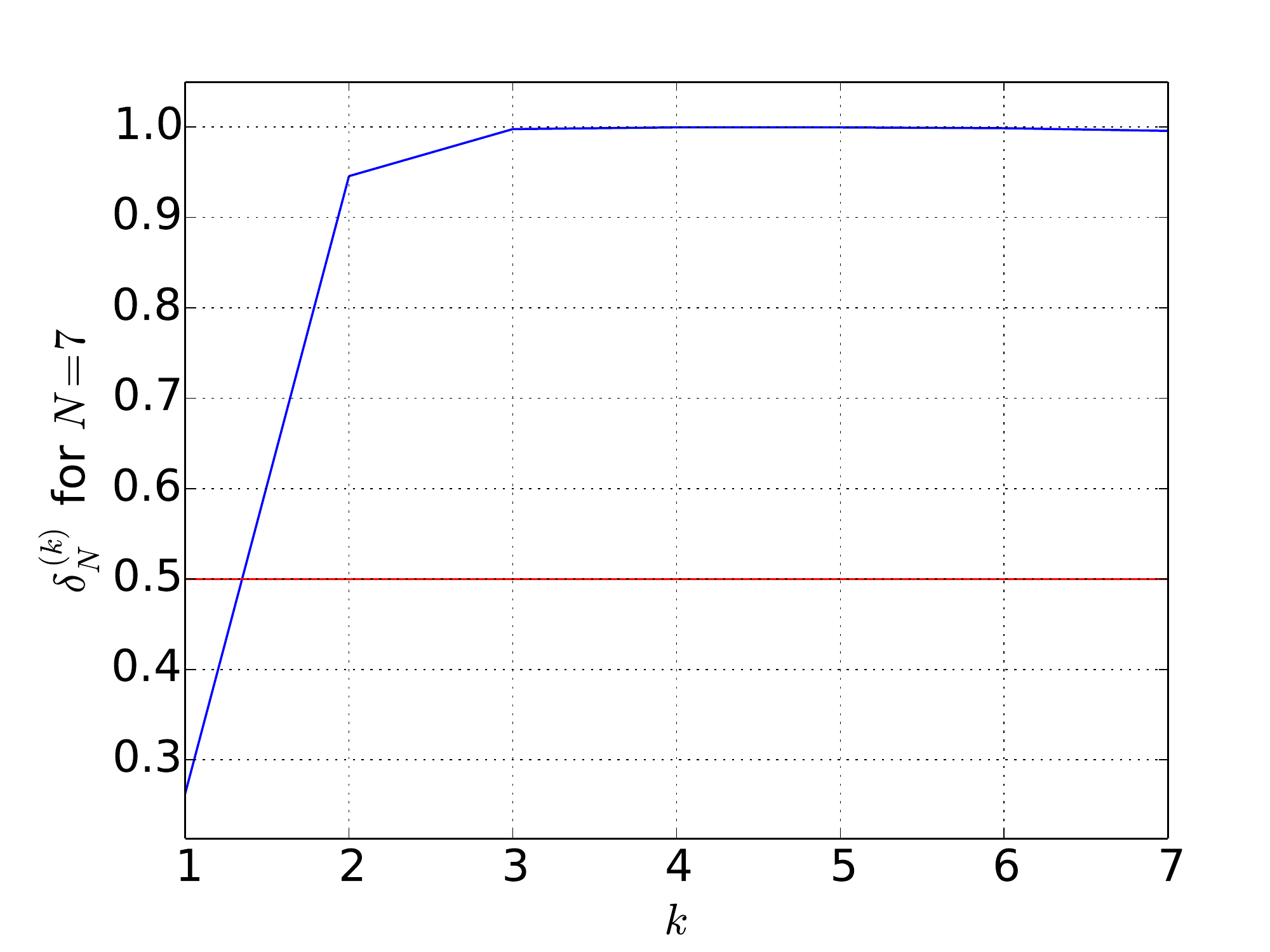}}
  \hspace{1pt}
  \subfloat[$N=10$]{\label{fig6b}\includegraphics[scale=0.29]{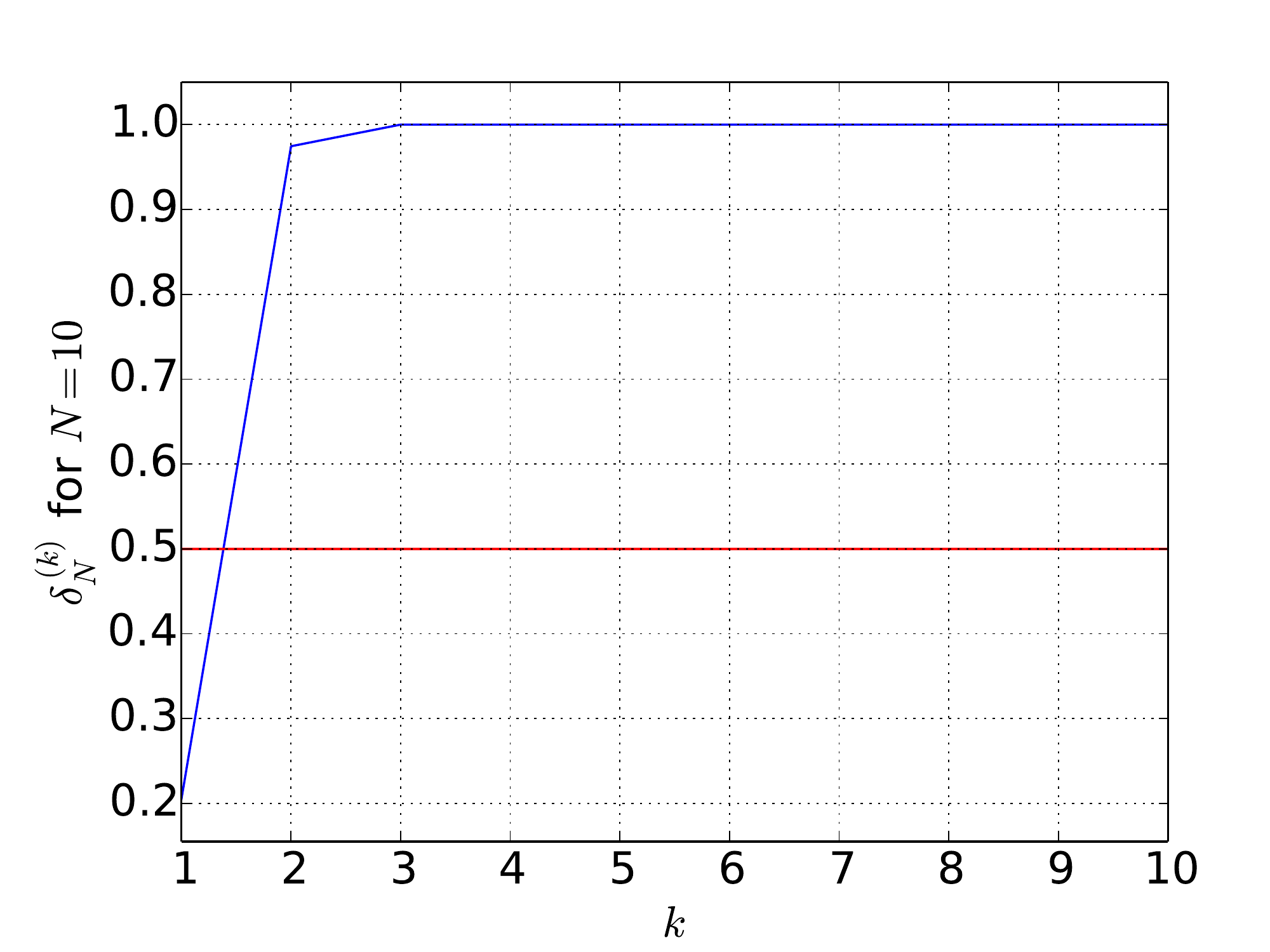}}
  \caption{The value of $\delta_N^{(k)}$ for $1\leq k\leq N$ and $N=7$, 10. The red line represents $1/2$. While for the basic mutational model, the process covers a little portion of the hypercube in $\mathcal O(N)$, which is smaller for bigger $N$, allowing more than one mutation at each step, the process can actually cover at least a half of the graph in a time of the same order.}
  \label{fig:delta}
\end{figure}

Then the strongest condition on $\delta$ is the one given by Proposition \ref{thm:epsexpa} (that we need to obtain  \eqref{eqtempbis})~: $\delta\leq1/2$. Therefore,  the 2-BRW-$\Pp^{(k)}$ is growing exponentially until it covers half of the hypercube. 
Then the way the rest of the hypercube is covered is not known.
\endproof

As we saw in the previous section, we are now able to prove an equivalent of Theorem \ref{lempartial} for this BRW:

\begin{theorem}\label{lempartialbis}
Given a simple 2-BRW-$\Pp^{(k)}$, there exists a time $T$ such that $T=\mathcal O(N)$ and with high probability $|S_T|\geq\delta2^{N}$, $\delta$ satisfying the hypothesis of Theorem \ref{lemnubis}.
\end{theorem}

In Figure \ref{fig:delta} we plot the value of the maximal proportion of vertices of the hypercube we can cover in $\mathcal O(N)$ considering a 2-BRW-$\Pp^{(k)}$. Of course, the case corresponding to $k=1$ ($\Pp^{(k)}=\Pp$) is obtained by Lemma \ref{lemnu}, and we denote $\delta_N^{(1)}:=2^{-r(N)}$ as obtained in  \eqref{eq:Rn}. These simulations shows that actually $\delta_N^{(k)}>1/2$ for all $k\geq2$ even for small $N$. This result suggests that once we break the bipartiteness by allowing at least two switch-type mutations at each time step, then the corresponding BRW invades at least half of the hypercube vertex set in $\mathcal O(N)$ (see Section \ref{sec:dis} for a further overlook on this issue). 

\begin{figure}[ht!]
  \centering
  \subfloat[$N=7$]{\label{fig6a}\includegraphics[scale=0.29]{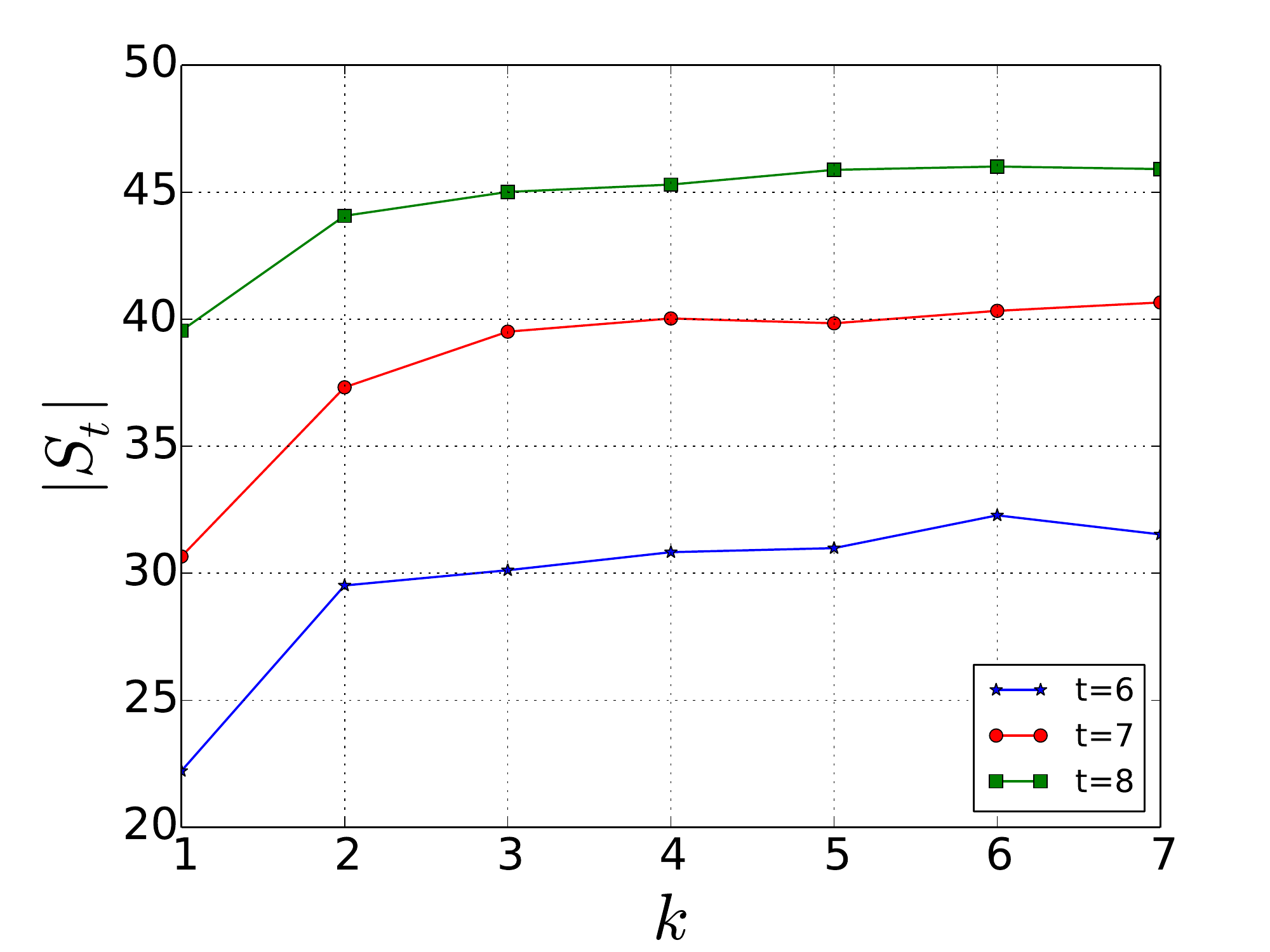}}
  \hspace{1pt}
  \subfloat[$N=10$]{\label{fig6b}\includegraphics[scale=0.29]{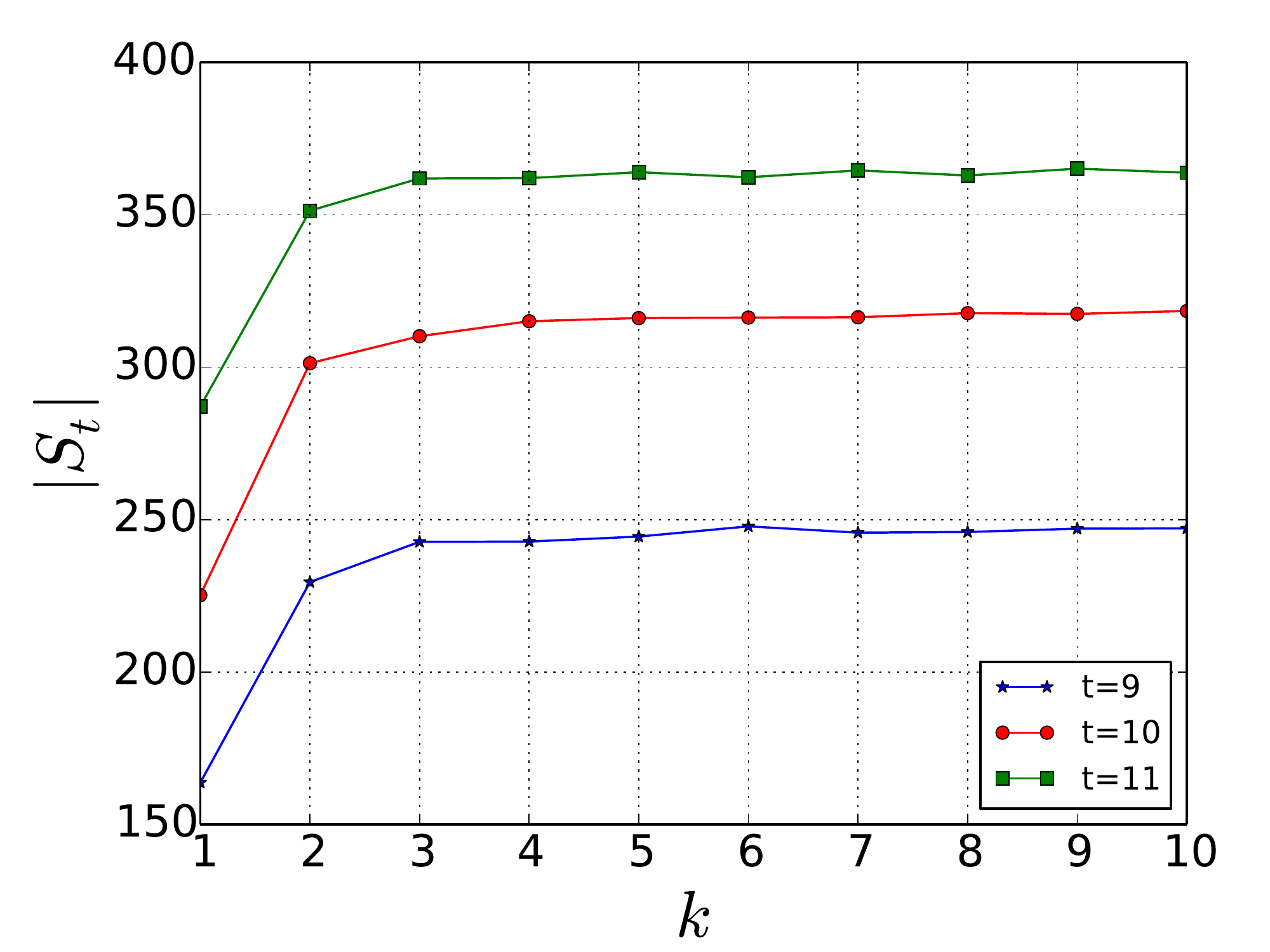}}
  \caption{Average size of $S_t$ after $t=N-1$, $t=N$ and $t=N+1$ time steps, comparing the 2-BRW-$\Pp^{(k)}$ with $k\,\in\,\{1,\dots,N\}$. Here we plot the average values obtained over 100 simulations.}
  \label{fig:PkaN}
\end{figure}

\begin{remark}
The definition of $\delta$  in Theorems \ref{lemnubis} and \ref{lempartialbis} does not depend on $k\geq2$ :  even for small values of $k$, we are able to cover at most a half of the hypercube vertex set in a time $\mathcal O(N)$. In Figure \ref{fig:PkaN} we simulated the average size of $S_t$ obtained by considering a 2-BRW-$\Pp^{(k)}$, with $k\,\in\,\{1,\dots,N\}$ for different time $t$. Simulations shows that the size of $S_t$ significantly increases passing from $\Pp$ to $\Pp^{(2)}$, and it is almost constant for $k$ between 3 and $N$. 
\end{remark}

\begin{remark}
The method applied here does not allow to prove a better covering of $\hyp$ in a time $T=\mathcal O(N)$ than the one obtained in Theorem \ref{lempartialbis} for matrix $\Pp^{(k)}$

\begin{lemma}\label{lem:maxT}
Let $\mathcal M$ be a transition probability matrix over $\hyp$, represents a $d$-regular, connected and non bipartite graph. 
Let $\lambda_2$ be the second largest eigenvalue of $\mathcal M$. 
Given a simple 2-BRW-$\mathcal M$, there exists a $\delta(\mathcal M):=\frac{de^{-2}+d-2}{(1-\lambda_2^2)(d^2e^{-2}+d-2)}-\frac{\lambda_2^2}{1-\lambda_2^2}$ such that in a time $T=\mathcal O(N)$ with high probability $|S_T|\geq\delta(\mathcal M)2^N$. For every such transition probability matrix  $\mathcal M$, $\delta(\mathcal M)\leq1/2$.
\end{lemma}

In other words, applying the method used in Sections \ref{sec3:2} and \ref{sec3:3}, the best result we can prove for a 2-BRW-$\mathcal M$  is $|S_T|\geq\delta2^N$, $\delta\leq1/2$  in a time $T=\mathcal O(N)$.

\proof
The assumptions made over $\mathcal M$ and Corollary \ref{cor:epsalphadelta} imply that $\mathcal M$ expresses a $(\alpha,\delta)$-expander graph, with $\delta\leq1/2$ and $\alpha=\left(\delta(1-\lambda_2^2)+\lambda_2^2\right)^{-1}$.
Let us consider a simple 2-BRW-$\mathcal M$. The method used in Sections \ref{sec3:2} and \ref{sec3:3} for $\Pp$ and $\Pp^{(k)}$ respectively allows to find $\delta(\mathcal M)$ (depending on $d$ and $\lambda_2$, as given in Lemma \ref{lem:maxT}) s.t. there exists a time $T=\mathcal O(N)$ s.t. with high probability $|S_T|\geq\delta(\mathcal M)2^N$. However we have a restriction over $\delta(\mathcal M)$, given by Proposition \ref{thm:epsexpa}, which is $\delta(\mathcal M)\leq1/2$.
\endproof

Furthermore, a similar threshold shall be explicit for a generic 2-BRW on a bipartite graph defined on the vertices of $\hyp$. 
At each time step we are observing the evolution of $|S_t|$ over a half part of $\hyp$, 
hence over a state-space of size $2^{N-1}$. Proceeding as above, we obtain the same results with $N-1$ instead of $N$. Therefore,  
the best result we can expect in a time $T=\mathcal O(N)$ is $|S_T|\geq\delta2^{N-1}$, $\delta\leq1/2$. 

\newcommand{\hypnmu}{\mathcal{H}_{N-1}}

\begin{figure}[hb!]
  \centering
  \includegraphics[scale=0.35]{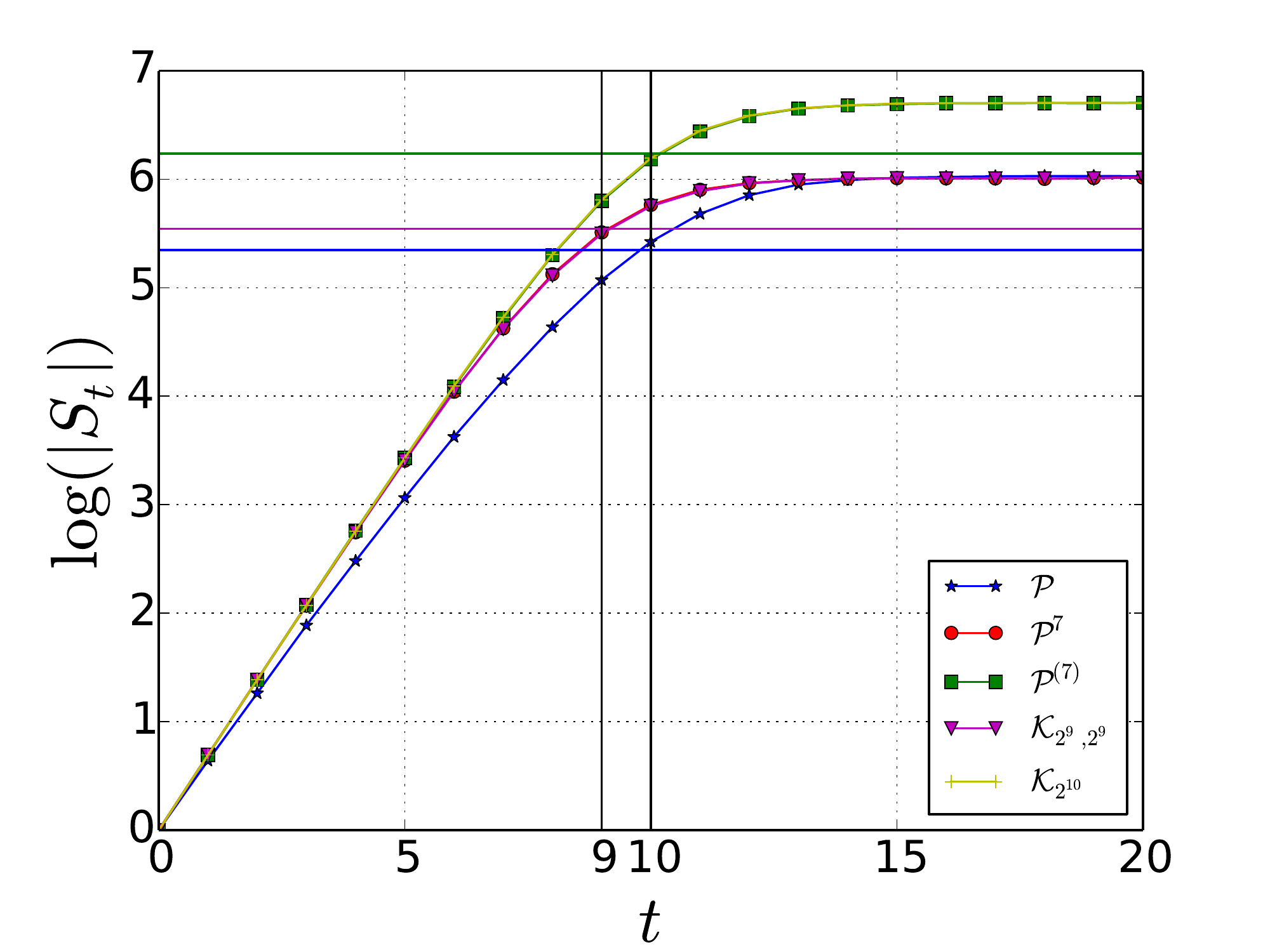}
  \caption{Evolution of $|S_t|$ on a log. scale. We compare the 2-BRW for different transition probability matrices and $N=10$. 
  Average values are plotted obtained over 40 simulations. The green horizontal line corresponds to $\log(2^9)$~: we can observe that in a time $t=N$ we do not overtake this threshold, even while considering the complet graph over $2^N$ vertices. 
  The magenta horizontal line corresponds to $\log(2^8)$~: BRWs associated to bipartite graphs do not cover more than this value in a time $t=N-1$. 
  Finally, the blue horizontal line represents the theoretical size of $S_t$ in a time $\mathcal O(N)$ for the simple 2-BRW-$\Pp$, as obtained in Theorem \ref{lempartial}. The curves corresponding to $\Pp^{(7)}$ and $\mathcal K_{2^{10}}$ are almost overlapping~:  the expansion properties of both matrices  ensure a covering of the same order. The same holds for curves corresponding to $\Pp^7$ and $\mathcal K_{2^9,2^9}$, which characterize bipartite graphs over $\{0,1\}^N$ with an appreciable vertex expansion.}
  \label{fig:St_exp}
\end{figure}

In Figure \ref{fig:St_exp} we test the evolution of the active set size for a 2-BRW corresponding to other transition probability matrices over $\hyp$ which assure good expansion properties. 
We show the ability of these simple 2-BRWs to cover $\hyp$ for $N=10$, in logarithmic scale. In particular we consider 5 transition probability matrices:
\begin{itemize}
\item $\Pp$, in blue.
\item $\Pp^{7}$, in red.
\item $\Pp^{(7)}$, in green.
\item $\mathcal K_{2^9,2^9}$, in magenta, defined as follows:
\begin{equation*}
\mathcal K_{2^9,2^9}:=\frac{1}{2^9}\left(\begin{array}{cc}
\boldsymbol{0}_{2^9} & \mathcal J_{2^9} \\
\mathcal J_{2^9} & \boldsymbol{0}_{2^9}\end{array}
\right),
\end{equation*}
where $\boldsymbol{0}_{2^9}$ is a $2^9\times2^9$ matrix with all entries 0 and $\mathcal J_{2^9}$ is a $2^9\times2^9$ matrix with all entries 1. This is the transition probability matrix corresponding to the complete bipartite graph on $2\times2^9$ vertices.
\item $\mathcal K_{2^{10}}:=\displaystyle\frac{1}{2^{10}-1}\left( \mathcal J_{2^{10}}-\mathcal I_{2^{10}}\right)$, in yellow, where $\mathcal I_{2^{10}}$ is the $2^{10}$-identity matrix. This transition probability matrix corresponds to the complete graph on $2^{10}$ vertices.
\end{itemize}

We  introduce  the complete bipartite graph and the complete graph in order to test the ability in invading the state-space for two transition probability matrices with strong expansion properties. This choice is not biologically motivated and we do not expect that they actually describe  actual mutation rules. \\

We observe that although for the complete graph, which has the best expansion property, 
in a time $t=N$ we can cover about a half of the state-space, as with the simple 2-BRW-$\Pp^{(k)}$. 
Even for small $t>0$, the process corresponding to $\Pp^{(7)}$ is faster when compared to 2-BRW-$\Pp^7$. It is interesting to compare this fact with a phenomenon observed in  \cite{balelli2015branching} where we investigated the typical time-scale of the exploration of $\hyp$ considering RWs without branching. We demonstrated that for $k>2$, $\Pp^{k}$ optimizes the hitting time to reach a certain configuration, if compared to $\Pp^{(k)}$. When we take into account the {\em branching} equivalent of these RWs, the exploration of $\hyp$ is more efficient using $\Pp^{(k)}$ as transition probability matrix instead of $\Pp^{k}$. That suggests that once  added a branching process, the oscillations due to bipartiteness are of greater amplitude
and forbid a quick covering even for small $t$. 
\end{remark}

\section{Extensions of the model}\label{sec5}

In this Section we set some variants of the model considered so far, in which we take into account the multiplicity of each vertex. 
This adds a further building block to our model.
 Indeed, taking into account the number of particles lying on the same vertex allows to consider the size of the effective population 
 and not only how many different BCR configurations are expressed at a certain time. 
 Moreover, considering multiplicity also allows us to have a better chance of making $|S_t|$ grow faster, 
 where $|S_t|$ represents here the number of vertices of $\{0,1\}^N$ on which at least one particle lies. In Section \ref{sec3:4} we consider BRWs with multiplicity and fixed number of offspring $c$ at each time step. 
 Then, in Section \ref{sec5:1}, we give to each individual a probability $p$ to divide~: 
 we observe the impact  of  division  on the limiting distribution. 
 Finally, in Section \ref{sec:dis}, we observe and discuss, 
 through  computer simulations, 
 a model for which the division rate depends on affinity.
 
\subsection{$\mathbf c$-BRW with multiplicity}\label{sec3:4}

 
At time $t\geq0$ we have exactly $c^t$ particles, as there is no death nor selection. 
We consider the distribution of these $c^t$ particles within $\hyp$. 
In order to do so, we define the Markov process $(X_t^i)_{t\geq0}$, where for all $i\,\in\,\{1,\dots,2^N\}$, 
$X_t^i$ corresponds to the number of particles lying on the $i^\textrm{th}$ node at time $t$. 
Proposition \ref{proPpar} is given in the more general case of a $c$-BRW with multiplicity on a given $d$-regular graph~: the case we are interested in is an application with $c=2$ and $d=N$.

\begin{proposition}\label{proPpar}
Given a $c$-BRW with multiplicity on a $d$-regular graph, then for all $s\geq0$~:
\begin{equation*}
\Pro\left[X_t^i=s\,\left\vert\displaystyle\sum_{j\sim i}X_{t-1}^j=n\right.\right]=\left\{\begin{array}{ll}
\displaystyle {c n\choose s} \frac{(d-1)^{c n-s}}{d^{c n}} & \textrm{if $s\leq c n$}, \\ \\
0 & \textrm{otherwise}.
\end{array}\right.
\end{equation*}
\end{proposition}

\proof
We show that conditioning on $\sum_{j\sim i}X_{t-1}^j=n$, $X_t^i$ follows a binomial distribution $\mathcal{B}\left(c n,\frac{1}{d}\right)$. For all $j\sim i$ let us define the random variables $\vec Z_{l,r}^j$, where $\vec Z_{l,r}^j$ corresponds to the vertex chosen by the $l^{\textrm{th}}$ particle lying on $j$ in its $r^{\textrm{th}}$ trial, with $ j\,\in\, S_{t-1}\cap \mathcal{N}(\{i\})$, $1\leq l\leq X_{t-1}^j$ and $1\leq r\leq c$. Then we have:
$$\Pro[\vec Z_{l,r}^j=i]=1/d\quad\forall\;j,l,r$$
At each trial of each particle lying on a vertex $j$, we have exactly a probability of $1/d$ of success (\emph{i.e.} going on vertex $i$) and a probability of $1-1/d$ of failure, and we have exactly $c n$ independent and identically distributed trials. Then the result follows.
\endproof

In particular, we have:

\begin{proposition}\label{proPpar3}
Given a $c$-BRW with multiplicity on the complete graph on $d$ vertices $\mathcal{K}_d$, the distribution of $X_t^i$ given $X_{t-1}^i=s'$ is a binomial distribution with parameters $c^t-cs'$ and $\frac{1}{d-1}$, \emph{i.e.} for all $s\geq0$~:
\begin{equation*}
\Pro[X_t^i=s\,|\,X_{t-1}^i=s']=\left\{\begin{array}{ll}
\displaystyle {c^t-cs'\choose s} \left(\frac{1}{d-1}\right)^s\left(1-\frac{1}{d-1}\right)^{c^t-cs'-s} & \textrm{if $s\leq c^t-cs'$}, \\ \\
0 & \textrm{otherwise}.
\end{array}\right.
\end{equation*}
\end{proposition}

Proposition \ref{proPpar3} shows that, for a complete graph on $N$ vertices, 
the probability of having $s$ particles at time $t$ on the $i^{\textrm{th}}$ 
node depends on the number of particles laying on $i$ at time $t-1$.

\proof
In this particular case, $i$ is connected to all nodes of the graph, except itself. Therefore each one of the $c^t$ particles produced at time $t$ has a probability $1/(d-1)$ to go to $i$~: we have to remove the particles that will leave from $i$, and this is exactly $cs'$.
\endproof

We establish another property of the $c$-BRW with multiplicity~: the asymptotic distribution of the $c^t$ individuals for $t\to\infty$.
This concludes this section. 

\begin{lemma}\label{lem_dis_lim_kBRW}
Let $\mathcal M$ be the transition probability matrix corresponding to a finite connected graph $G=(V,E)$, $\vec m$ its stationary distribution. Let us suppose $\mathcal M$ aperiodic, and let us consider a $c$-BRW-$\mathcal M$ starting from a generic initial distribution $\vec p$. Therefore:
\begin{equation*}
\forall\,i\,\in\,V,\; \frac{X_t^i}{c^t}\to\vec m_i \textrm{ in probability, for $t\to\infty$.}
\end{equation*}
\end{lemma}

\proof
The position of each  of the $c^t$ individuals at time $t$ corresponds to the position reached by a RW with $\mathcal M$ as transition probability matrix, starting from the initial distribution $\vec p$ and independently form  others individuals. 
In other words, at time $t$ we are considering the position of $c^t$ parallel RWs-$\mathcal M$ starting from the same initial distribution. 
For all $j\,\in\,\{1,\dots,c^t\}$, let $(X_{j,t})_{t\geq0}$ i.i.d RWs with transition probability matrix $\mathcal M$ and starting from the initial distribution $\vec p$. By hypothesis, for all $i\,\in\,V$, $\Pro(X_{j,t}=i)\to\vec m_i$ for $t\to\infty$. The result follows since convergence in law to a constant implies convergence in probability.
\endproof

\begin{remark}
Numerically, we compare the average size of $S_t$ for $t=N=10$ for the simple 2-BRW-$\Pp$, the simple 2-BRW-$\mathcal K_{2^9,2^9}$ and the 2-BRW-$\Pp$ with multiplicity. Table \ref{tab:hist2} below shows the average values obtained over 100 simulations. As expected, the 2-BRW-$\Pp$ with multiplicity is faster than the simple 2-BRW-$\Pp$ because of the number of particles within the population, which is not affected nor by selection or death, neither by coalescence. At each step, each particle can divide and colonize a new vertex of the hypercube, therefore we have a better chance to cover faster a half of the state-space (we recall that $\Pp$ is a bipartite graph). Moreover, we can observe that the simple 2-BRW-$\mathcal K_{2^9,2^9}$ is faster than the simple 2-BRW-$\Pp$~: indeed $\mathcal K_{2^9,2^9}$ has better expander properties, and thus the BRW invades more efficiently the state-space as noticed in Sec. \ref{sec3}. 

\begin{table}[ht!]
\caption{Average size of $S_t$ after 10 time steps, comparing the simple 2-BRW-$\Pp$, the simple 2-BRW-$\mathcal K_{2^9,2^9}$ and the 2-BRW-$\Pp$ with multiplicity. We denote by $\widehat{|S_{10}|}_n$ the average value obtained over $n$ simulations and by $\widehat{\sigma}_n$ its corresponding estimated standard deviation.}\label{tab:hist2}
\centering
\begin{tabular}{l c c c c}
\hline\noalign{\smallskip}
\textbf{Model} & $\boldsymbol{N}$ & $\boldsymbol{n}$ & $\boldsymbol{\widehat{|S_{10}|}_n}$ & $\boldsymbol{\frac{\widehat{\sigma}_n}{\sqrt{n}}}$ \\
\noalign{\smallskip}\hline\noalign{\smallskip} 
\textbf{Simple 2-BRW-$\mathbf \Pp$} & 10 & 100 & 222.36 & 3.376\\
\textbf{Simple 2-BRW-$\mathbf{\mathcal K_{2^9,2^9}}$} & 10 & 100 & 318.04 & 1.231\\
\textbf{2-BRW-$\mathbf \Pp$ with multiplicity} & 10 & 100 & 398.42 & 0.972\\
\noalign{\smallskip}\hline
\end{tabular}
\end{table}

\end{remark}

\subsection{Limiting distribution for the BRW-$\Pp$ with multiplicity and division rate $p$.}\label{sec5:1}

Lemma \ref{lem_dis_lim_kBRW} can not be applied to the 2-BRW-$\Pp$ with multiplicity. Indeed, the bipartite structure of the corresponding graph prevents the convergence through the stationary distribution, \emph{i.e.} the homogeneous probability distribution, which we denote by $\boldsymbol \pi$ \cite{balelli2015branching}. We can overcome this problem by considering a BRW-$\Pp$ with multiplicity and with a non constant division rate $p$. 

\begin{definition}\label{def:proba:p}
Let us fix $p\,\in\,]0,1[$. The process starts with a single individual located on an arbitrary node of $\hyp$. Each time step, a particle lying on a certain node $\vec x_i$ of $\hyp$ gives rise to $2$ daughter cells and die with probability $p$. 
With probability $1-p$, it remains in the population for the next time step. 
When division occurs, each newborn particle choses a neighbor node according to matrix $\Pp$, independently and with replacement, and move on it.
\end{definition}

The introduction of a division rate has two immediate consequences. First, it slows down population's growth. 
In order to evaluate the expected number of individuals at time $t$, 
we consider a generic Galton-Watson process (\cite{harris2002theory}, chapter I).

\begin{proposition}\label{prop:esp_GW}
Let $Z_t$ be the r.v. describing the number of individuals at generation $t$ starting from $Z_0=1$ individual. 
We assume that each individual divides indepently from the others and from previous generations.
Let $\vec p:=(p_k, k=0,1,2,\dots)$ be a probability distribution s.t. 
$p_k$ gives the probability of having $k$ offsprings in the next generation. 
At each time step, given $Z_t=k$, $Z_{t+1}$ behaves as $k$ independent copies of $Z_1$. 
Therefore: $\E(Z_t)=\left(\E(Z_1)\right)^t$.
\end{proposition}

In our specific case we have:
\begin{itemize}
\item $p_1=1-p$
\item $p_2=p$
\item $p_k=0$ for all $k\neq1,2$
\end{itemize}

Which gives:

\begin{equation}\label{eq:EZt}
\E(Z_t)=(1+p)^t<2^t\textrm{ as $p<1$.}
\end{equation}

\begin{remark}
One can observe that $Z_t=\displaystyle\sum_{i=1}^{2^N}X_t^i$, where $X_t^i$ describes the number of individuals lying on vertex $i$ at time $t$.
\end{remark}

The addition of the parameter $p$ overcomes issues related to the bipartite structure of the graph, 
discussed in Section \ref{sec:bip}.

\begin{lemma}\label{lem:bip_p}
Let us consider a $BRW$ with multiplicity on a finite connected bipartite graph $G_b$. Let $\vec p:=(p_k, k=0,1,2,\dots)$ be the probability distribution of the number of offsprings of each individuals for the next generation, s.t. $p_1>0$ and $p_0+p_1<1$. Then there exists a time $t\geq0$ and two nodes $\vec x_1$, $\vec x_2$ s.t. $\vec x_1\sim\vec x_2$ and $\vec x_1$, $\vec x_2\,\in\,S_t$.
\end{lemma}

Lemma \ref{lem:bip_p} implies that for this type of BRWs, independently from the bipartite structure of $G_b=(V,E)$, 
there exists a time $t>0$ s.t. $S_t=V$ (see Section \ref{ss:BRW_nonbip}).

\proof
Let $0<T<\infty$ s.t. $\vec x_i\,\in\,S_T$ ($T$ exists as $G_b$ is finite and connected). As $p_0+p_1<1$, $\exists\,k\geq2$ s.t. $p_k>0$. Then with probability $p_k$, $\exists\,\vec x_{i,1},\dots,\vec x_{i,k}\,\in\,\mathcal{N}(\{\vec x_i\})$ s.t. $\{\vec x_{i,1},\dots, \vec x_{i,k}\}\,\in\,S_{T+1}$. As $p_1>0$, with positive probability at least one among these $k$ vertices does not divide: let $\overline k\,\in\{1,\dots,k\}$ s.t. $\vec x_{i,\overline k}\,\in\,S_{T+2}$. Moreover w.p.p. one among $\{\vec x_{i,1},\dots, \vec x_{i,k}\}\setminus\{\vec x_{i,\overline k}\}$ divides and w.p.p. one of its offsprings migrates to $\vec x_{i}$. Therefore, w.p.p. $\{\vec x_{i,\overline k}, \vec x_i\}\,\in\,S_{T+2}$, and $\vec x_{i,\overline k}\sim \vec x_i$.
\endproof

\noindent We give an equivalent of Lemma \ref{lem_dis_lim_kBRW} for BRWs characterized by Definition \ref{def:proba:p}.

\begin{lemma}\label{lem_dis_lim_kBRW2}
Let $\mathcal M$ be the transition probability matrix corresponding to a finite connected graph $G=(V,E)$, $\vec m$ its stationary distribution. Let us consider a BRW-$\mathcal M$ with multiplicity starting from a generic initial distribution. Let $\vec p:=(p_k, k=0,1,2,\dots)$ be the probability distribution of the number of offsprings of each individual for the next generation, with $p_1>0$ and $p_0+p_1<1$. We denote by $Z_t$ the r.v. describing the population size at generation $t$ (starting from $Z_0=1$). For all $i\,\in\,V$ let $X_t^i$ be the r.v. describing the number of individuals lying on vertex $i$ at time $t$. Therefore:
\begin{equation*}
\forall\,i\,\in\,V,\; \frac{X_t^i}{\left(\E(Z_1)\right)^t}\to\vec m_i \textrm{ in probability for $t\to\infty$.}
\end{equation*}
\end{lemma}

\proof
The proof is the same as for Lemma \ref{lem_dis_lim_kBRW}. In this case, we do not need the hypothesis of aperiodicity of $\mathcal M$ as the problem of an eventual periodicity is overcome by the addition of the distribution of the number of offspring $\vec p$, as shown in Lemma \ref{lem:bip_p}. 
\endproof

Lemma \ref{lem_dis_lim_kBRW2} allows us to prove:

\begin{corollary}\label{cor_dis_lim_kBRW2}
Let us consider a BRW-$\Pp$ with multiplicity and division rate $p\,\in\,]0,1[$.
\begin{equation*}
\forall\,i\,\in\,\{0,1\}^N,\; \frac{X_t^i}{(1+p)^t}\to\frac{1}{2^N} \textrm{ in probability for $t\to\infty$.}
\end{equation*}
\end{corollary}

\proof
We have already determined $\E(Z_t)$ corresponding to the BRW-$\Pp$ with multiplicity and division rate $p$ (cf. \eqref{eq:EZt}). Therefore, in order to prove Corollary \ref{cor_dis_lim_kBRW2} we have just to observe that the stationary distribution for $\Pp$ is the homogeneous probability distribution on $\{0,1\}^N$. Then the result follows  applying Lemma \ref{lem_dis_lim_kBRW2}.
\endproof

\begin{remark}\label{rem:PL}
In a previous paper \cite{balelli2015branching} we overcame the problem of the bipartiteness of the graph underlined by $\Pp$ by adding $N$ loops at each node. That corresponds to take into account matrix $\Pp_L:=\frac{1}{2}(\Pp+\mathcal I_{2^N})$ instead of $\Pp$. Considering a BRW-$\Pp$ with multiplicity and division rate $p=1/2$ is equivalent to consider a 2-BRW-$\Pp_L$ with multiplicity, but with coalescence of those offsprings which decide to remain in place. The only difference is the size of the population at time $t$, which is $2^t$ in the case of a 2-BRW-$\Pp$ with multiplicity and $(3/2)^t$ in the other case. The choice of $\Pp_L$ as transition probability matrix has also biological motivations. Indeed division of B-cells in GCs is asymmetric \cite{meyer2012theory,barnett2012asymmetric}: only one between the two daughter cells has a mutated trait. 
\end{remark}

\begin{remark}
More  generally, let us consider a transition probability matrix $\mathcal M$ on a graph $G=(V,E)$, with $|V|=n$. We can see a BRW-$\mathcal M$ with multiplicity and division rate $p$ as a 2-BRW-$\mathcal M_p$ with multiplicity, where $\mathcal M_p:=p\mathcal M+(1-p)\mathcal I_{n}$. Of course, we need to take the same caution as in Remark \ref{rem:PL} about the number of individuals at time $t$.
\end{remark}

\subsection{BRW-$\Pp$ with multiplicity and affinity dependent division}\label{sec:dis}

In previous sections, the limiting distribution of traits (with or without division rate) only depends on the stationary distribution of the considered transition probability matrix. 
In particular, if the stationary distribution is homogeneous, 
than for $t$ big enough all individuals are uniformly distributed over the state-space. 
From a biological point of view, it does not seem so efficient to explore all the state-space. 
It will be rather more interesting to drive mutations through the region of the state-space with greater affinity for the target trait.
We can therefore propose a model in which we introduce a division rate dependent on the affinity of the cell. \\

\begin{figure}[ht!]
\centering
\subfloat[ ]{\includegraphics[width = 3in]{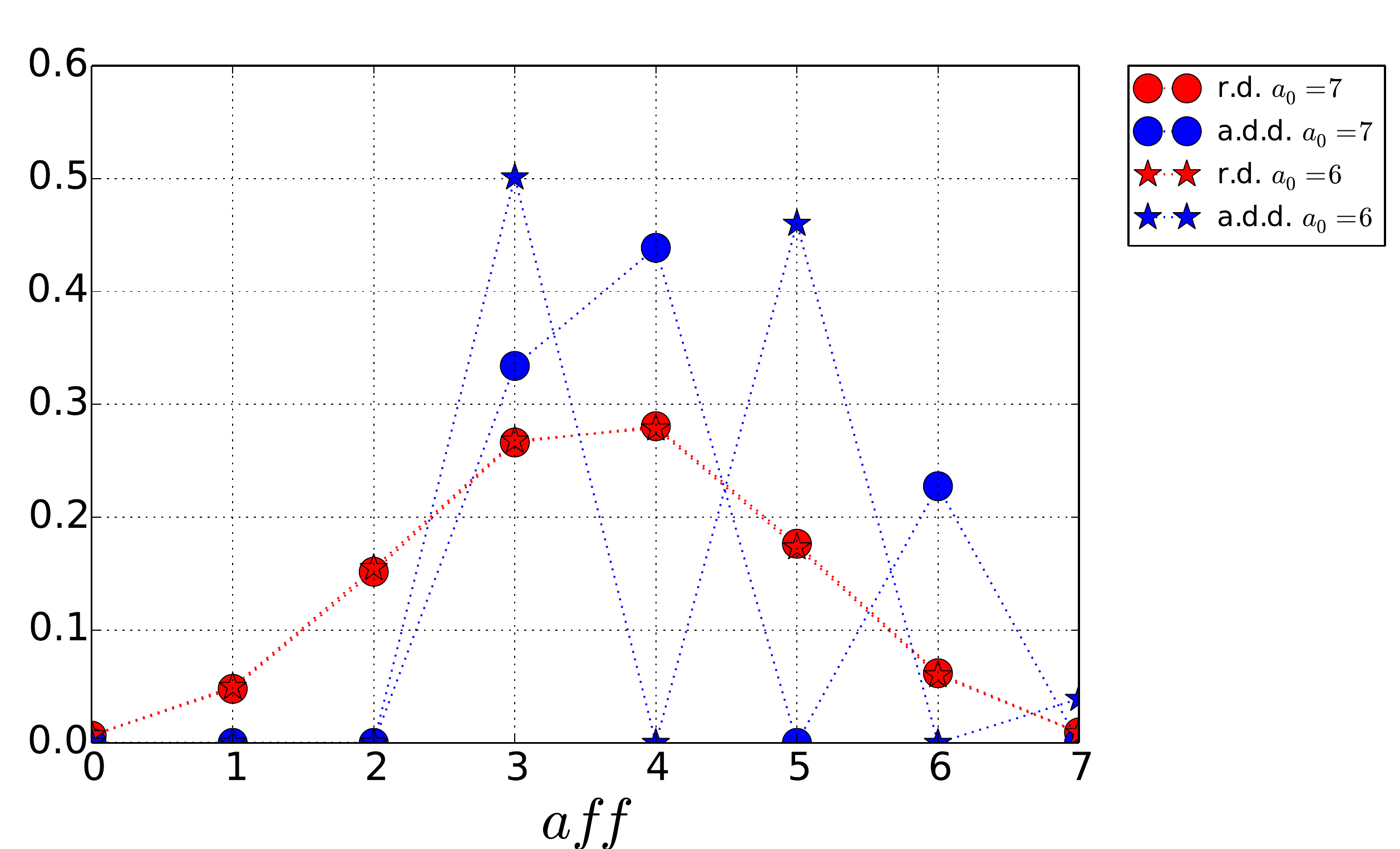}\label{fig1:intro}}\\
\subfloat[ ]{\includegraphics[width = 2in]{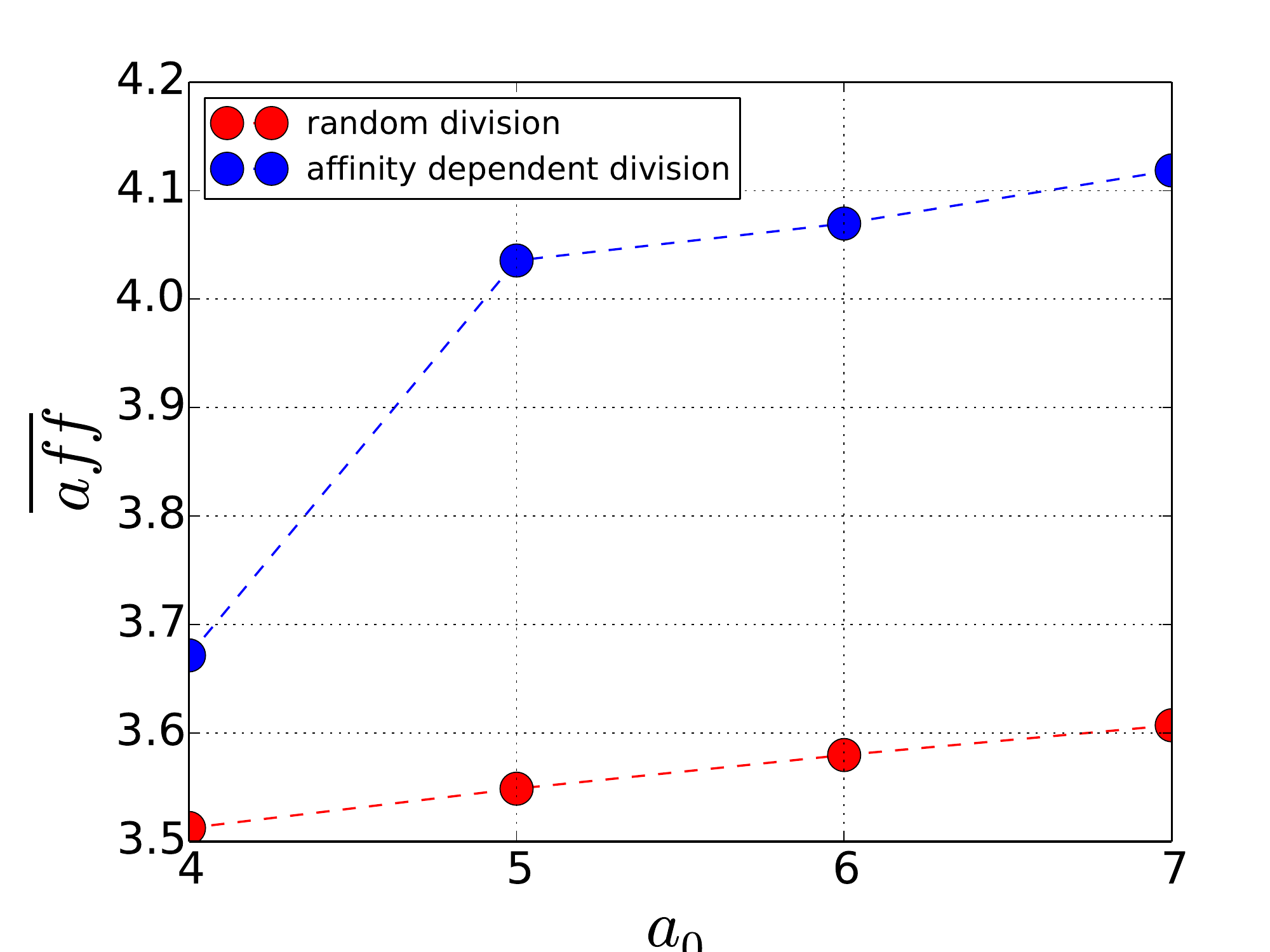}\label{fig2:intro}}
       ~ 
\subfloat[ ]{\includegraphics[width = 2in]{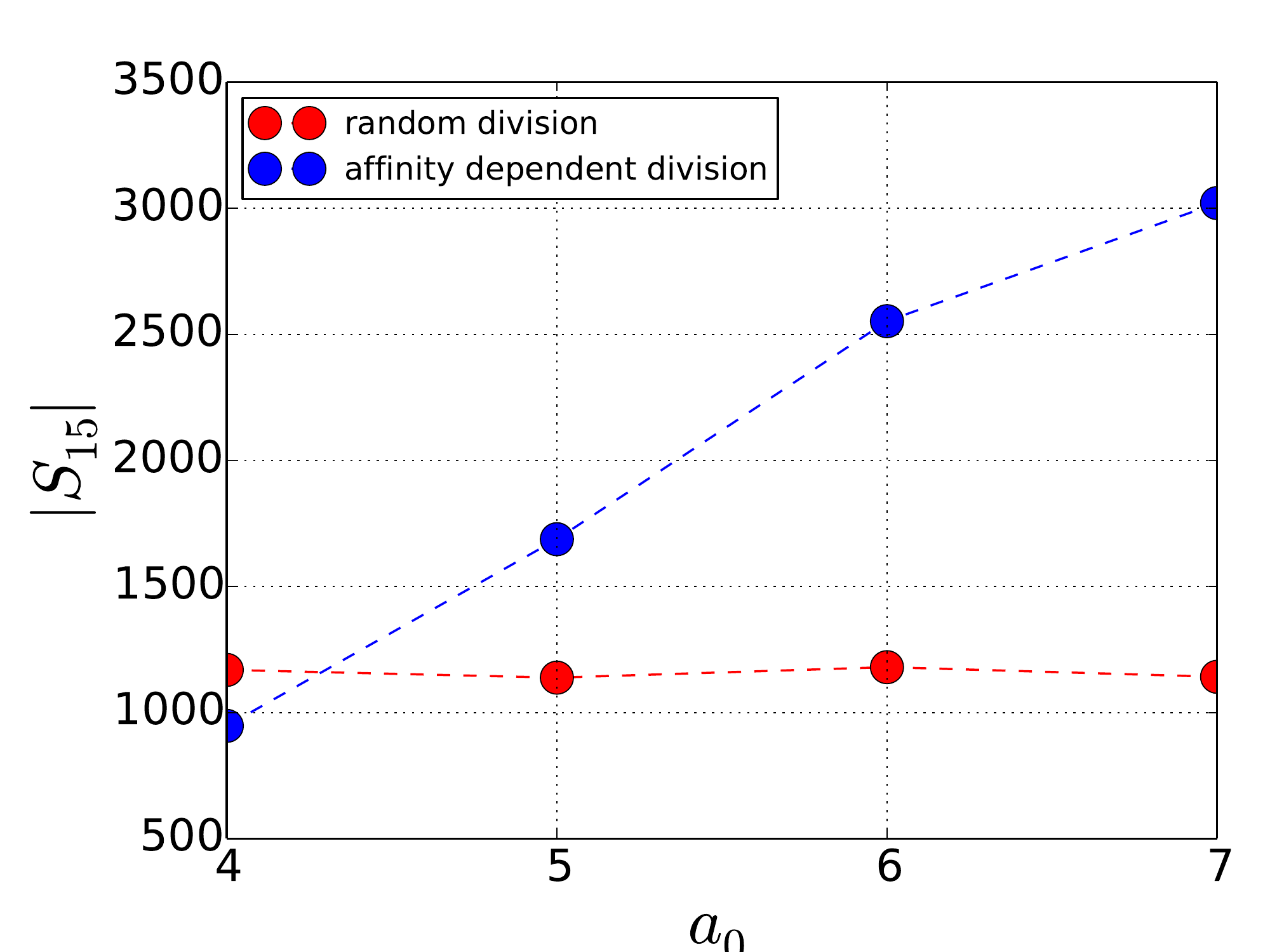}\label{fig2:intro}}
\caption{Simulations of the BRW-$\Pp$ with multiplicity, comparing a model with division rate $p=0.6$ (in red) and a model with affinity dependent division (in blue). In this last case individuals having affinity at least 4 with the target vertex divide and mutate accordingly to matrix $\Pp$, they remain unchanged in the population otherwise. (a) Distribution of the affinity to the antigen after 15 time steps, starting from initial affinity 7 (circles) and 6 (stars) respectively. (b) Dependance of the average affinity (after 15 time steps) on the initial affinity $a_0$. (c) Dependance of the final population size (after 15 time steps)  on the initial affinity $a_0$.}\label{fig:test_mult_comp}
\end{figure}

Formally, $\forall\,\vec x_i\,\in\,\hyp$, let $p_d(\vec x_i)$ be the probability of division of an individual lying on vertex $\vec x_i$. We can define an increasing function $f$ s.t. $p_d(\vec x_i)=f\left(\textrm{aff}(\vec x_i,\cible)\right)$, where $\textrm{aff}(\vec x_i,\cible)=N-h(\vec x_i,\cible)$ is the affinity of $\vec x_i$ with respect to the target trait $\cible$ (Definition \ref{def:affhamming}), and $h$ return the Hamming distance. The aim is to be able to privilege those individuals having better fitness. This choice has biological motivations. Indeed, recent evidence shows that during GC reaction the acquisition of highest affinity for the presented antigen regulates proliferation and diversification of B-cells \cite{gitlin2014clonal}. 
In Figure \ref{fig:test_mult_comp} we compare a model of BRW-$\Pp$ with multiplicity and division rate $p=0.6$ with a model of BRW-$\Pp$ with multiplicity and affinity dependent division. In this case, we chose a very simple function for the division rate, defined $\forall\,\vec x_i\,\in\,\hyp$, as follows:
\begin{equation}
p_d(\vec x_i)=\left\{\begin{array}{ll}
0 & \textrm{if $\textrm{aff}(\vec x_i,\cible)<N-\overline h_s$} \\
1 & \textrm{if $\textrm{aff}(\vec x_i,\cible)\geq N-\overline h_s$}
\end{array}\right.
\end{equation}
We plot results obtained for $N=7$ and $\overline h_s=3$~: all individuals having affinity at least 4 with the target trait divide and mutate accordingly to matrix $\Pp$, they remain unchanged in the population otherwise.\\

In Figure \ref{fig:test_mult_comp} (a) we represent the final distribution of the affinity of the traits within the population after 15 time steps. 
As expected, the distribution corresponding to the first model is binomial and does not depend on the initial Hamming distance. 
Indeed, from Corollary \ref{cor_dis_lim_kBRW2} we know that the distribution of the traits is uniform on $\{0,1\}^N$. 
We have just to remark that in $\{0,1\}^N$ there are exactly ${N \choose h}$ nodes having Hamming distance $h$ from a given vertex, $0\leq h\leq N$~: this determines the proportion of individuals having a given affinity after 15 time steps.
The support of the  distribution at time step 15 for the second model corresponds to vertices 
having affinity 3, 4 or 6 (resp. 3, 5, 7) with the target trait for an initial affinity $a_0=7$, (resp. $a_0=6$). 
Indeed, as $a_0\geq 4$, the total population can be divided in two subpopulations. 
The sub-population  whose affinity with the target trait is greater than 4 follows a standard 2-BRW-$\Pp$ with multiplicity. 
Therefore, we can observe the effects of the bipartiteness of the graph: only traits whose affinity has the same parity as $a_0$ are expressed at even time step. On the contrary, at odd time steps only vertices with affinity having the opposite parity as $a_0$ are expressed.
The other sub-population is composed by those individuals that after an unfavorable mutation obtain a trait having affinity exactly 3. They remain unchanged for all further time steps, as they can not divide nor die. Therefore, through further time steps, individuals with affinity 3 can only continue to accumulate. This is due to the definition of $p_d(\vec x_i)$ as a step function.\\

Figure \ref{fig:test_mult_comp} (b) shows the average affinity of the population after 15 time steps. We can see that for the BRW-$\Pp$ with division rate $0.6$ this depends very lightly from the initial affinity, while, as expected, the initial affinity strongly influences the final one if we allow only individuals having affinity greater than 3 to divide. Finally in Figure \ref{fig:test_mult_comp} (c) we see the size of the population after 15 time steps. Again, in the case of random division with rate 0.6, the initial affinity does not affect the final population size, which is always approximately $1.6^{15}\simeq1152.92$.

\section{Conclusions and perspectives}\label{sec:conclusion}

In this paper, we introduce and study BRWs on binary strings, modeling the evolution of cells in a mutation-division process. The edge set (or graph) associated to $\hyp:=\{0,1\}^N$, hence the corresponding transition probability matrix, reflects  mutations allowed during the evolutionary process. Graph's characteristics determine the behavior of the BRW, \emph{e.g.} its ability in covering $\hyp$ or the limiting distribution of the traits, 
as shown in Sections \ref{sec3} and \ref{sec5}.\\

We particularly focus on the expander property of the graphs when giving quantitative results about the expected portion of $\hyp$ covered in $\mathcal O(N)$. 
We observe that strong expansion properties enable a faster invasion of the state-space. 
From a biological point of view, this property is significant since it ensures that starting from one or a few B-cell, the GC can produce, hence test a huge variety of BCRs against the target antigen. 
Indeed, GCs seem to be oligoclonal \cite{kroese1987germinal,maclennan1994germinal}, 
which means that they develop from very few initial naive B-cells (three, on average). 
Therefore, starting from a single clonal population, it is of interest  to understand how a B-cells population  invades the BCR state-space. \\

For this reason, in Section \ref{sec3}, we consider the state-space $\hyp$ of every possible $N$-length string (modeling B-cell traits), and compare the ability of different mutation rules in colonizing $\hyp$ in a time $\mathcal O(N)$. We develop upon a method used in \cite{duttacoalescing} to evaluate partial cover times on expander graphs. Nevertheless, our approach differs from \cite{duttacoalescing}. Indeed, we fix the state-space
and  the main question becomes :
how many nodes we are able to activate in a time $\mathcal O(N)$ for a given graph?
In particular, we observe that while matrix $\Pp$, which denotes the structure of the standard $N$-dimensional hypercube, can cover a quite small portion of $\hyp$ in a time $\mathcal O(N)$, 
the mutation rule $\Pp^{(k)}=\frac{1}{k}\sum_{i=1}^k\Pp^i$ leads to a significantly bigger expansion which does not strongly depends on 
$k$, for values of $k$ greater than 2. \\

In Section \ref{sec3}, we  show that if we simply consider the expansion properties of the structure built over $\hyp$,  the covering in $\mathcal O(N)$ is limited at a half the state-space (Lemma \ref{lem:maxT}). 
This favors the hypothesis that the expansion property is not enough to insure a quick covering of a large portion of the state-space~: 
considering self-avoiding BRWs on connected graphs could be  more efficient, although these are not necessarily good expanders. 
On the other hand, from a biological point of view, 
it may not be so efficient to explore the whole state-space, but rather to steer mutations toward a specific region of the state-space with the best affinity. 
Indeed, the production of new clones has a cost in terms of time and energy, therefore it does not make sense to produce a huge variety of cells with any possible fitness with the presented antigen. Models considered in this paper share this drawback~: even if a bigger portion of possible traits is expressed in a time $\mathcal O(N)$, we can not say much about their average fitness.\\

We can propose many possible solutions to this problem. We can for example privilege individuals with good fitness by considering a model with affinity dependent division, as discussed in Section \ref{sec:dis}. Another possibility is to consider transition probability matrices whose stationary distribution is concentrated on a specific region of the state-space containing the fittest traits. Indeed, as we observe in Section \ref{sec3:4}, given this hypothesis than the distribution of traits for a 2-BRW with multiplicity only depends on the stationary distribution of the  transition probability matrix under consideration. 
In this case the problem is~: does this matrix accounts for realistic mutations? 
Another way to drive mutations towards a specific region of the state-space is, 
of course, the introduction of a selection mechanism, which we study in a separate work \cite{balelli2016GW}.



\bibliographystyle{apt}      

\end{document}